\documentclass[11pt]{article}
\usepackage{graphicx}
\usepackage{epsfig}
\usepackage{latexsym,amsfonts,amsbsy,amssymb}
\usepackage{amsmath,amsthm}
\usepackage{hyperref}
\usepackage{color}
\usepackage{float} 
\newcommand{\calF}{{\cal F}}
\newcommand{\bR}{\mathbb R}
\newcommand{\R}{\mathbb R}
\newcommand{\bC}{\mathbb C}
\newcommand{\Z}{\mathbb Z}
\newcommand{\balpha}{\boldsymbol{\alpha}}
\newcommand{\bbeta}{\boldsymbol{\beta}}
\begin{document}
\title{Taylor-Fourier approximation}
\author{
M.P. Calvo\footnote{IMUVA and Departamento de Matem\'atica Aplicada, Facultad de Ciencias, Universidad de
Valladolid, Spain. E-mail: mariapaz.calvo@uva.es}\,, J. Makazaga\footnote{Universidad del Pais Vasco (UPV/EHU), Donostia-San Sebasti\'an, Spain. E-mail: joseba.makazaga@ehu.eus}\, and A. Murua\footnote{Universidad del Pais Vasco (UPV/EHU), Donostia-San Sebastián, Spain. E-mail: ander.murua@ehu.eus}}

\maketitle

\begin{abstract}
In this paper, we introduce an algorithm that provides approximate solutions to semi-linear ordinary differential equations with highly oscillatory solutions, which, after an appropriate change of variables, can be rewritten as non-autonomous systems with a $(2\pi/\omega)$-periodic dependence on $t$. The proposed approximate solutions are given in closed form as functions $X(\omega t,t)$, where $X(\theta,t)$ is (i) a truncated Fourier series in $\theta$ for fixed $t$ and (ii) a truncated Taylor series in $t$ for fixed $\theta$, which motivates the name of the method.

These approximations are uniformly accurate in $\omega$, meaning that their accuracy does not degrade as $\omega \to \infty$. In addition, Taylor-Fourier approximations enable the computation of high-order averaging equations for the original semi-linear system, as well as related maps that are particularly useful in the highly oscillatory regime (i.e., for sufficiently large $\omega$).

The main goal of this paper is to develop an efficient procedure for computing such approximations by combining truncated power series arithmetic with the Fast Fourier Transform (FFT). We present numerical experiments that illustrate the effectiveness of the proposed method, including applications to the nonlinear Schrödinger equation with non-smooth initial data and a perturbed Kepler problem from satellite orbit dynamics.
\end{abstract}
\section{Introduction}
\label{sec:introduction}
In this paper  we are interested in systems of ordinary differential equations of the form
\begin{equation}
\label{eq:TF0}
\frac{d}{dt}{\bf x} = \mathcal{A} {\bf x} + {\bf g}({\bf x}),
\end{equation}
where ${\bf g}: \bR^D \to \bR^D$ is a smooth map and $\mathcal{A}$ is a real $D\times D$ matrix with purely imaginary eigenvalues.

Although the core concepts and techniques can be naturally extended to the general case of arbitrary purely imaginary eigenvalues, 
we primarily focus on the fully resonant case, in which the ratio of any pair of eigenvalues is rational. In this scenario, the solution ${\bf x}(t) = e^{t \mathcal{A}} {\bf x}(0)$ to (\ref{eq:TF0}) with ${\bf g} \equiv {\bf 0}$ becomes $(2\pi/\omega)$-periodic in $t$, where $\omega$ is the fundamental frequency (i.e., all eigenvalues of $\mathcal{A}$ are integer multiples of $\omega i$.)

We thus assume that $\mathcal{A}=\omega \, A$, where $A$ is a real $D\times D$ matrix whose eigenvalues are integer multiples of the imaginary unit  $i$.
We focus on approximating the solution ${\bf x}(t)$ of

\begin{equation}
\label{eq:TF1_1}
\frac{d}{dt}{\bf x} = \omega A {\bf x} + {\bf g}({\bf x}),
\end{equation}
 supplemented with the initial condition
\begin{equation}
\label{eq:TF_icond}
{\bf x}(0)={\bf x}_0 \in \bR^D.
\end{equation}
 Examples of differential systems \eqref{eq:TF1_1} with these properties can be found, for instance,  when dealing with semi-discretized dispersive partial differential equations~\cite{C-GP_2015,H-O2010,Trefethen}, or in some formulations of the equations of motion in satellite dynamics~\cite{Stiefel-Scheifele}.

In the present work, we aim at obtaining a closed form approximation to the solution  ${\bf x}(t)$ of (\ref{eq:TF1_1}) that exhibits uniform accuracy across varying frequencies $\omega$ for a time-interval (centered at $t=0$) of length independent of $\omega$.

With the change of variables
\begin{equation}
\label{eq:change_var}
{\bf x}(t) = e^{t \omega A} {\bf y}(t),
\end{equation}
the initial value problem (\ref{eq:TF1_1})--(\ref{eq:TF_icond}) reads
\begin{equation}
\label{eq:TF1_2}
\frac{d}{dt}{\bf y} = {\bf f}(\omega t, {\bf y}), \quad {\bf y}(0)={\bf y}_0,
\end{equation}
with
\begin{equation}
\label{eq:TF1_3}
{\bf f}(\theta, {\bf y}) = e^{-\theta A} {\bf g}(e^{\theta A} {\bf y}),
\end{equation}
which is $(2\pi)$-periodic in $\theta$, and ${\bf y}_0 = {\bf x}_0$. 

In the context of problems of the form \eqref{eq:TF1_1} obtained from the space semi-discretization of dispersive PDEs, the numerical integration with a Runge-Kutta method of the transformed system \eqref{eq:TF1_2}, as a means to solve (\ref{eq:TF1_1})--(\ref{eq:TF_icond}), leads to the so-called Lawson's generalized Runge-Kutta methods~\cite{Lawson,H-O2010}.

A completely different context where problems of the considered form arise is the numerical orbit propagation of artificial satellites.
The equations of motion of the satellite can be brought into the form (\ref{eq:TF1_1}), for instance, by using the Kustaanheimo-Stiefel (KS) formulation~\cite{K-S1965}. In this context,  (\ref{eq:TF1_2}) is referred to as the VOP (variation of parameters) formulation of the system (\ref{eq:TF1_1}). In particular, a VOP formulation of the KS equations was  presented and analyzed by Stiefel and Scheifele in \cite{Stiefel-Scheifele}.  Other VOP formulations (obtained by using different sets of coordinates) of the equations of motion in satellite dynamics are also of the form \eqref{eq:TF1_2}.

 The Taylor-Fourier approximations that we propose here draw inspiration from previous works on the analysis of numerical integration of oscillatory problems:

 \begin{itemize}
\item  The technique of modulated Fourier expansions~ (see e.g. \cite{C-H-L,H-L-W},  and references therein), when applied to (\ref{eq:TF1_2}) would lead to a formal expansion
 \begin{equation}
 \label{eq:MFE}
\sum_{k \in \Z} e^{i k \omega t} {\bf y}_k(t)
\end{equation}
of the solution ${\bf y}(t)$ of (\ref{eq:TF1_2}), where ${\bf y}_k(t)$ ($k\in \Z$)  are smoothly varying functions (i.e., their derivatives are bounded independently of $\omega$). An approximation of ${\bf y}(t)$ can be obtained from (\ref{eq:MFE}) by replacing each ${\bf y}_k(t)$ by its $d$th degree truncated Taylor expansion centered at $t=0$, and by discarding the terms corresponding to Fourier modes with $|k|>M$, that is,
\begin{equation}
 \label{yd}
{\bf y}(t) \approx \sum_{k=-M}^{M} e^{i k \omega t} \sum_{j=0}^d  t^j {\bf y}_{k,j}.
\end{equation}
The $(M,d)$-Taylor-Fourier approximations that we propose in the present paper are precisely of the form (\ref{yd}).

\item  In~\cite{Ch-M-SS2012}, the high-order averaging of systems of the form (\ref{eq:TF1_2}) is studied by expanding its solution ${\bf y}(t)$ as a formal series indexed by multi-indices $(k_1 \cdots k_{\ell}) \in \Z^{\ell}$. Such series expansion can be best described by parametrizing the equations (\ref{eq:TF1_2}) with  a dummy variable $\epsilon$ as
\begin{equation}
\label{eq:TF1_2_epsilon}
\frac{d}{dt}{\bf y} = \epsilon \, {\bf f}(\omega t, {\bf y}), \quad {\bf y}(0)={\bf y}_0.
\end{equation}
  Its solution ${\bf y}(t)$ is expanded as
\begin{equation*}
{\bf y}(t) = {\bf y}_0 + \sum_{\ell=1}^{\infty} \epsilon^{\ell} \sum_{(k_1\cdots k_{\ell}) \in \Z^{\ell}} \alpha_{k_1\cdots k_{\ell}}(t) \, {\bf f}_{k_1 \cdots k_{\ell}}({\bf y}_0),
\end{equation*}
where each ${\bf f}_{k_1 \cdots k_{\ell}}$ is a map that is independent of $\omega$, and each $\alpha_{k_1 \cdots k_{\ell}}(t)$ is a scalar function that only depends on $t$ and $\omega$ and the multi-index
$(k_1, \ldots, k_{\ell})$, but is independent from the actual map ${\bf f}(\theta,{\bf y})$. More precisely, each $\alpha_{k_1 \cdots k_{\ell}}(t)$ is a linear combination with rational coefficients of terms of the form $\omega^{j-\ell} t^j e^{i \, k \, \omega\, t}$, $0\leq j \leq \ell$. An equivalent expansion with series indexed by mode-colored rooted trees was derived in~\cite{Ch-M-SS2012a}.
After collecting terms with equal powers of $\epsilon$ such expansions of the solution ${\bf y}(t)$ can be written as
\begin{equation}
\label{eq:y_epsilon}
{\bf y}(t) = {\bf y}_0 + \sum_{\ell=1}^{\infty} \epsilon^{\ell}
\sum_{k=-\infty}^{\infty} e^{i k \omega t} \sum_{j=0}^{\ell}  t^j {\bf y}_{k,\ell,j}.
\end{equation}
As shown in~\cite[Theorem 3.1]{Ch-M-SS2012}, under suitable conditions on the map ${\bf f}(\theta, {\bf y})$,  there exists $\eta>0$ independent of $\omega$ such that the series (\ref{eq:y_epsilon}) is absolutely and uniformly convergent for $|\epsilon\, t| \leq \eta$. An essential insight that inspired the consideration of Taylor-Fourier approximations, as presented in this paper, is the realization~\cite{Ch-M-SS2012a} that the approximations derived from the expansion (\ref{eq:y_epsilon}) by truncating terms involving higher powers of $\epsilon$,
\begin{equation}
\label{eq:y_epsilon_N}
{\bf y}^{[d]}(t) = {\bf y}_0 + \sum_{\ell=1}^{d} \epsilon^{\ell} \sum_{k=-\infty}^{\infty} e^{i k \omega t} \sum_{j=0}^{\ell} t^j {\bf y}_{k,\ell,j},
\end{equation}
can be systematically constructed by combining Picard iterations with the truncation of the series expansion in powers of $\epsilon$. Indeed, starting with ${\bf y}^{[0]}(t)=  {\bf y}_0$, for $d=1,2,3,\ldots$, one first determines
\begin{equation}
\label{eq:iter_1}
{\bf z}^{[d]}(t) =  \sum_{\ell=0}^{d} \epsilon^{\ell} \sum_{k=-\infty}^{\infty} e^{i k \omega t} \sum_{j=0}^{\ell} t^j {\bf z}_{k,\ell,j}
\end{equation}
uniquely such that
\begin{equation}
\label{eq:iter_2}
{\bf f}(\omega t, {\bf y}^{[d]}(t)) = {\bf z}^{[d]}(t)+ \mathcal{O}(\epsilon^{d+1}),
\end{equation}
and then obtains
\begin{equation}
\label{eq:iter_3}
{\bf y}^{[d+1]}(t) = {\bf y}_0 + \epsilon \int_{0}^{t}{\bf z}^{[d]}(\tau) d\tau.
\end{equation}

In the particular case where $\epsilon = 1$, one obtains a series expansion of the solution ${\bf y}(t)$ of (\ref{eq:TF1_2}) that is absolutely and uniformly convergent for $|t| \leq \eta$, where $\eta > 0$ is independent of $\omega$.

Taylor-Fourier approximations, as presented in this paper, emerged as a practical computational method for approximating such expansions. To construct approximations of the form (\ref{yd}), an iterative process similar to (\ref{eq:iter_1})--(\ref{eq:iter_3}) is employed. However, instead of truncating the series in powers of $\epsilon$, truncation is applied to powers of $t$, and the exact Fourier series expansion is replaced with trigonometric interpolation.
   \end{itemize}
Our Taylor-Fourier approximations are also related to the methods proposed in~\cite{B-C-M-R2017} for a more restricted class of problems of the form (\ref{eq:TF1_1}).


The aim of the present paper is: (i) to introduce and give a precise definition of Taylor-Fourier approximations, (ii) to relate them to high order stroboscopic averaging, (iii) to propose a procedure to compute them efficiently by making use of truncated power series arithmetic and the FFT algorithm, and (iv) to illustrate the application of the proposed approximation method  to certain well known problems of interest.

For given $M$ and $d$, the algorithm introduced in this paper provides closed-form expressions
\begin{equation}
\label{xd}
{\bf x}(t) \approx e^{t \omega A} \left(\sum_{k=-M}^M e^{ik\omega t} \sum_{j=0}^{d} t^j \, {\bf y}_{k,j}^{[d]} \right)
\end{equation}
of approximate solutions to (\ref{eq:TF1_1})-(\ref{eq:TF_icond}). These approximations maintain uniform accuracy in $\omega$ and remain valid over a time interval independent of $\omega$. Once computed, these approximations can be efficiently evaluated at arbitrary times $t$ without additional costly computations. This contrasts with standard numerical integration methods, which produce approximations to the exact solution at a finite sequence of intermediate times $0 < t_1 < t_2 < \cdots \leq T$, and whose accuracy typically deteriorates for large values of $\omega$.

The Taylor-Fourier approximations presented in this work are not intended to compete with numerical integration methods in terms of computational efficiency. Rather, their primary purpose is to provide reliable reference solutions over large time intervals. Additionally, they can serve as a valuable tool for assessing the accuracy and effectiveness of methods that approximate the solutions of high-order averaged equations.

 We have carried out a careful implementation of the Taylor-Fourier methods in Julia language, paying special attention to computational aspects such as memory management and efficient truncated power series  arithmetic.

In Section~\ref{sec:def} we give a precise definition of $(M,d)$-Taylor-Fourier approximations.
Section \ref{sec:implementation} is devoted to describe an algorithm designed to efficiently compute truncated Taylor-Fourier series.
 In Section~\ref{sec:averaging}, we explore the relationship of the Taylor-Fourier technique with high-order stroboscopic averaging.  We demonstrate in Section~5 the applicability of Taylor-Fourier approximations to problem classes extending beyond real semi-linear problems of the form \eqref{eq:TF0} with fully resonant frequencies.
In Section \ref{sec:experiments}  we illustrate the application of the Taylor-Fourier algorithm to a couple of relevant examples.   Finally, in Section~\ref{sec:conclusions}, we summarize our conclusions.

\section{Taylor-Fourier approximations}
\label{sec:def}

To lay the groundwork for the definition of fully discrete $(M,d)$-Taylor-Fourier approximations (\ref{yd}) to the solution ${\bf y}(t)$ of (\ref{eq:TF1_2}) for arbitrary pairs of positive integers $(M,d)$, we first introduce the semi-discrete $(\infty,d)$-Taylor-Fourier approximations, which correspond to the limit case as $M \to \infty$.

\subsection{Semi-discrete Taylor-Fourier approximations}
\label{subsec:semi-discrete}
For each positive integer $d$, the $(\infty,d)$-Taylor-Fourier approximation to the solution ${\bf y}(t)$ of (\ref{eq:TF1_2}) is of the form
\begin{equation}
\label{eq:Yd}
{\bf Y}^{[d]}(\omega t, t) =  \sum_{j=0}^{d} t^j \, {\bf Y}_j^{[d]}(\omega t),
\end{equation}
where each component of the $D$-vector functions ${\bf Y}_j^{[d]}(\theta)$ is $(2\pi)$-periodic in the angular variable $\theta$.
The function
 \begin{equation*}
{\bf Y}^{[d]}(\theta,t) =  \sum_{j=0}^{d} t^j \, {\bf Y}_j^{[d]}(\theta)
\end{equation*}
 is defined recursively. Starting with ${\bf Y}^{[0]}(\theta,t) = {\bf y}_0$, ${\bf Y}^{[d+1]}(\theta,t)$ is determined for $d\geq 0$ as follows:
\begin{enumerate}
\item First, we determine the $d$th-order truncated expansion in powers of $t$ for ${\bf f}(\theta, {\bf Y}^{[d]}(\theta,t))$,
\begin{equation}
\label{eq:Iter_1}
{\bf Z}^{[d]}(\theta,t) =  \sum_{j=0}^{d} t^j {\bf Z}_j^{[d]}(\theta) =
{\bf f}\left(\theta,  \sum_{j=0}^{d} t^j \, {\bf Y}_j^{[d]}(\theta)
\right) + \mathcal{O}(t^{d+1}).
\end{equation}
\item Second, we determine ${\bf Y}^{[d]}(\theta,t)$ from the identity
\begin{equation}
\label{eq:Iter_2}
{\bf Y}^{[d]}(\omega t, t) = {\bf y}_0 + \int_{0}^{t}{\bf Z}^{[d]}(\omega \tau, \tau) d\tau.
\end{equation}
\end{enumerate}

Notice that if, for each $j$, the function ${\bf Y}_j^{[d]}(\theta)$ admits a Fourier expansion
\begin{equation*}
{\bf Y}_j^{[d]}(\theta) = \sum_{k=-\infty}^{\infty} e^{i k \theta} {\bf y}^{[d]}_{k,j},
\end{equation*}
then the $(\infty,d)$-Taylor-Fourier approximation
\begin{equation}
\label{eq:Yd2}
{\bf Y}^{[d]}(\omega t,t) =  \sum_{k=-\infty}^{\infty} e^{i k \omega t} \sum_{j=0}^{d} t^j \, {\bf y}_{k,j}^{[d]}
\end{equation}
takes the form (\ref{yd}) in the limiting case $M\to\infty$.

\subsection{Fully discrete Taylor-Fourier approximations}
\label{subsec:fully_discrete}

For the time being, let us assume that the Fourier expansion of ${\bf f}(\theta,{\bf y})$ contains only a finite number of nonzero modes. That is, there exists a positive integer $M_0$ such that
\begin{equation}
\label{eq:finiteness}
{\bf f}(\theta,{\bf y}) = \sum_{k=-M_0}^{M_0} e^{i k \theta}\, \hat{\bf f}_k({\bf y}).
\end{equation}
In particular, if each component of ${\bf g}({\bf x})$ in \eqref{eq:TF1_1} is a polynomial function of ${\bf x} \in \mathbb{R}^D$, then (\ref{eq:TF1_3}) admits such a finite Fourier expansion.

If condition (\ref{eq:finiteness}) holds, then the semi-discrete Taylor-Fourier approximation (\ref{eq:Yd2}) can also be expressed as a finite sum. More precisely, it can be shown that ${\bf y}^{[d]}_{k,j}=0$ for all $|k| > M = d M_0$.

In this case, the first step (\ref{eq:Iter_1}) in the recursive definition of $(\infty,d)$-Taylor-Fourier approximations
can be split into two sub-steps. In the first sub-step,  truncated expansions in powers of $t$ of univariate functions of $t$ are performed, while in the second one, we apply trigonometric interpolation for uniformly spaced nodes in $[0,2\pi]$.
More precisely, we can consider $\theta_n = 2n \pi/(2M+1)$, $n=0, 1, 2, \ldots, 2M$, with  $M\geq (d+1) M_0$, and
\begin{itemize}
\item first compute for each node $\theta_n$ the vectors  ${\bf Z}^{[d]}_j(\theta_n)={\bf Z}_{n,j}^{[d]} \in \mathbb{C}^D$ from the truncated series expansion in powers of $t$
\begin{equation}
\label{eq:Iter_1_M}
 {\bf f} \left (\theta_n, \sum_{j=0}^{d} t^j \, {\bf Y}_j^{[d]}(\theta_n) \right ) = \sum_{j=0}^{d} t^j \, {\bf Z}_{n,j}^{[d]} + \mathcal{O}(t^{d+1})
\end{equation}
\item and then, for each $j=0,1,2,\ldots,d$, uniquely determine the functions
\begin{equation}
\label{eq:Iter_2_M}
{\bf Z}_j^{[d]}(\theta) = \sum_{k=-M}^{M} e^{i k \theta} {\bf z}_{k,j}^{[d]}
\end{equation}
from the conditions
\begin{equation*}
{\bf Z}_j^{[d]}(\theta_n)  = {\bf Z}_{n,j}^{[d]}, \quad \theta_n = \frac{2n\pi}{2M+1}, \quad n=0, 1, 2, \ldots, 2M.
\end{equation*}
\end{itemize}
The expansion (\ref{eq:Iter_1_M}) can be computed by evaluating the left-hand side of  (\ref{eq:Iter_1_M}) with truncated power series arithmetic. As for the trigonometric interpolation, it is efficiently performed with the application of the Discrete Fourier Transform (DFT) (see Section~\ref{sec:implementation}). Since the efficiency of DFT-based trigonometric interpolation is optimal for number of nodes that are powers of 2, interpolation with an even number of nodes is preferred in practice. We thus replace the interpolatory conditions (having an odd number of nodes) above by the alternative conditions with an even number of nodes,
\begin{equation}
\label{eq:Z_interp}
{\bf Z}_j^{[d]}(\theta_n)  = {\bf Z}_{n,j}^{[d]}, \quad \theta_n = \frac{n\pi}{M}, \quad n=0, 1, 2, \ldots, 2M-1,
\end{equation}
together with
\begin{equation}
\label{eq:Iter_4_M}
{\bf z}_{-M,j}^{[d]} = {\bf z}_{M,j}^{[d]}, \quad j=0,1,\ldots,d.
\end{equation}
Condition \eqref{eq:Iter_4_M} is imposed to ensure uniqueness of the trigonometric interpolation functions \eqref{eq:Iter_2_M}.

To sum up, \eqref{eq:Iter_1_M}--\eqref{eq:Iter_4_M} is a practical way of performing the first step \eqref{eq:Iter_1} in the recursive definition of $(\infty,d)$-Taylor-Fourier approximations to the solution ${\bf y}(t)$  for (\ref{eq:TF1_2}) provided that the finiteness condition \eqref{eq:finiteness} and $M>(d+1)M_0$ hold. Otherwise, the functions ${\bf Z}^{[d]}_j(\theta)$ determined by \eqref{eq:Iter_1_M}--\eqref{eq:Iter_4_M} do not satisfy \eqref{eq:Iter_1} exactly, as they will be affected by interpolation errors, that will depend on the actual value of the positive integer $M$.
In that general scenario, the function
\begin{equation}
\label{eq:Yd_M}
{\bf y}^{[M,d]}(t) = {\bf Y}^{[M,d]}(\omega t, t) =  \sum_{j=0}^{d} t^j \, {\bf Y}_j^{[M,d]}(\omega t),
\end{equation}
defined recursively by \eqref{eq:Iter_1_M}--\eqref{eq:Iter_4_M} and \eqref{eq:Iter_2} (with the superscript $[d]$ replaced by $[M,d]$) can be seen as a fully discrete substitute of the $(\infty,d)$-Taylor-Fourier approximation. We will say that \eqref{eq:Yd_M} is the $(M,d)$-Taylor-Fourier approximation to the solution ${\bf y}(t)$  for (\ref{eq:TF1_2}).

\section{Efficient computation of Taylor-Fourier approximations}
\label{sec:implementation}

Next, we describe an efficient method to compute the new coefficients ${\bf y}_{k,j}^{[M,d+1]} \in \mathbb{C}^D$ for $0 \leq j \leq d+1$ and $-M \leq k \leq M$, given the coefficients ${\bf y}_{k,j}^{[M,d]} \in \mathbb{C}^D$ for $0 \leq j \leq d$ and $-M \leq k \leq M$. For simplicity, we omit the superscripts $[M,d]$ in the remainder of this section. This computation proceeds in four sequential steps.

\subsection{Step 1: Evaluation of trigonometric polynomials}

For $n=0, 1, \ldots, 2M-1$, we calculate the vector coefficients
\begin{equation}
\label{eq:Ynj}
{\bf Y}_{n,j} =  {\bf Y}_{j}(\theta_n) = \sum_{k=-M}^M e^{i k n \pi/M}{\bf y}_{k,j},
\quad j=0, 1, \ldots, d,
\end{equation}
(with the nodes $\theta_n$ considered in \eqref{eq:Z_interp})
by applying the inverse DFT
$${\bf Y}_{n,j} = \frac{1}{2M}\sum_{k=0}^{2M-1} e^{ikn\pi/M} \hat{\bf Y}_{k,j}, \quad n=0, 1, \ldots, 2M-1, $$
where for $j=0,1,\ldots,d$,
$$
\begin{array}{l}
\hat{\bf Y}_{0,j} = 2M \, {\bf y}_{0,j}, \quad \hat{\bf Y}_{M,j} = 2M ({\bf y}_{M,j} + {\bf y}_{-M,j}), \\[8pt]
\hat{\bf Y}_{k,j} = 2M \, {\bf y}_{k,j}, \quad \hat{\bf Y}_{2M-k,j} = 2M{\bf y}_{-k,j}, \quad k=1, \ldots, M-1.
\end{array}
$$
\subsection{Step 2: Truncated Taylor expansions}

For each  $n=0, 1, \ldots, 2M-1$,
we compute ${\bf Z}_{n,j}$, $j=0, 1, \ldots, d$, such that
\begin{equation}
\label{eq:f_expansion}
{\bf f} \left ( {\theta_n, }\sum_{j=0}^{d} t^j \, {\bf Y}_{n,j} \right ) = \sum_{j=0}^{d} t^j \, {\bf Z}_{n,j} + {\cal O}(t^{d+1}).
\end{equation}

 The required truncated Taylor expansion \eqref{eq:f_expansion} can be computed efficiently without resorting to symbolic calculations by evaluating the univariate function on the left-hand side of \eqref{eq:f_expansion} using truncated power series arithmetic. To achieve this, one can implement the mapping ${\bf f}(\theta,{\bf y})$ as a function that takes an angle $\theta$ and $D$ polynomials in $t$ of degree $d$, returning the corresponding $D$ polynomials of degree $d$. This approach has been used in both examples presented in Section~\ref{sec:experiments}. Alternatively, numerous libraries are available across different computing environments that support arithmetic for truncated power series.

\subsection{Step 3:  Trigonometric interpolation}
For each $j=0, 1, \ldots ,d$, we determine the coefficients ${\bf z}_{k,j}$, $-M \leq k \leq M$, of the trigonometric polynomial
$${\bf Z}_j(\theta) = \sum_{k=-M}^M e^{ik\theta}{\bf z}_{k,j},$$
that satisfy the interpolatory conditions
$${\bf Z}_j(\theta_n) = {\bf Z}_{n,j}, \quad n=0, 1, 2, \ldots, 2M-1,$$
together with ${\bf z}_{-M,j} = {\bf z}_{M,j}$.

This task can be efficiently done by computing the discrete Fourier Transform (DFT)
\begin{equation}
\label{eq:dft}
\hat{\bf Z}_{k,j} = \sum_{n=0}^{2M-1} e^{-ikn\pi/M} \, {\bf Z}_{n,j}, \quad k=0, 1, \ldots, 2M-1,
\end{equation}
and setting
\begin{equation}
\label{hat_z}
\begin{array}{l}
{\bf z}_{0,j} = \dfrac{1}{2M} \hat{\bf Z}_{0,j}, \quad {\bf z}_{M,j} = {\bf z}_{-M,j} = \dfrac{1}{4M} \hat{\bf Z}_{M,j}, \\[8pt]
{\bf z}_{k,j} = \dfrac{1}{2M} \hat{\bf Z}_{k,j}, \quad {\bf z}_{-k,j} = \dfrac{1}{2M} \hat{\bf Z}_{2M-k,j}, \quad k=1, \ldots, M-1.
\end{array}
\end{equation}

\noindent{\em Remark:} For systems of the form \eqref{eq:TF1_2} in the real domain, the symmetries present in the complex vector resulting from the application of DFT to a real vector can be used to compute it more efficiently (with the real FFT algorithm) than in the general case of complex vectors. In particular, this implies that each ${\bf z}_{-k,j}$ coincides with the complex conjugate of ${\bf z}_{k,j}$, which allows us to halve the number of coefficients that  {must} be computed and saved.

\subsection{Step 4: Computing quadratures}

Given
$${\bf z}(t)= \sum_{k=-M}^M e^{ik\omega t} \left ( \sum_{j=0}^{d} t^j \, {\bf z}_{k,j} \right ),$$
we want to compute
\begin{equation}
\label{integral}
\begin{split}
{\bf y}(t) &= {\bf y}_0 + \int_0^t {\bf z}(s) ds\\
&=t^{d+1} \, {\bf y}_{0,d+1} + \sum_{k=-M}^M e^{ik\omega t} \left ( \sum_{j=0}^{d} t^j \, {\bf y}_{k,j}\right ).
\end{split}
\end{equation}

First, we notice that ${\bf z}(t)$ can be written as
\begin{equation}
\label{z_p}
{\bf z}(t) = {\bf p}(t) + \sum_{j=0}^{d} t^j \, \tilde{{\bf Z}}_{j}(\omega t),
\end{equation}
where
\begin{equation}
\label{p_and_hat_Z}
{\bf p}(t) = \sum_{j=0}^{d} t^j \, {\bf z}_{0,j}, \quad
\tilde{{\bf Z}}_{j}(\theta) = \sum_{k=1}^M (e^{ik\theta} {\bf z}_{k,j} + e^{-ik\theta} {\bf z}_{-k,j}), \quad j=0, \ldots, d.
\end{equation}
Therefore,
\begin{equation}
\label{y_q}
{\bf y}(t) = {\bf y}_0 + {\bf q}(t) + \sum_{j=0}^d t^j \, \tilde{{\bf Y}}_{j}(\omega t),
\end{equation}
with
\begin{eqnarray}
\label{hat_Y}
\tilde{{\bf Y}}_{j}(\theta) &=& \sum_{k=1}^M (e^{ik\theta} {\bf y}_{k,j} + e^{-ik\theta} {\bf y}_{-k,j}), \quad j=0, \ldots, d, \\
\label{poly_q}
{\bf q}(t) &=& - \tilde{{\bf Y}}_{0}(0) + \int_0^t {\bf p}(s) ds.
\end{eqnarray}
The coefficients ${\bf y}_{k,j}, {\bf y}_{-k,j}$, $j=0, 1, \ldots,  {d+1}$, $k=1, \ldots, M$,
can be efficiently obtained as follows:
\begin{itemize}
\item Set $\tilde{{\bf Y}}_{d+1}(\theta) \equiv {\bf 0}$, and notice that taking derivatives in \eqref{y_q} with respect to $\theta=\omega t$, and using \eqref{integral} and \eqref{z_p} one gets
$$\omega \frac{d}{d\theta}\tilde{{\bf Y}}_{j}(\theta) = \tilde{{\bf Z}}_{j}(\theta) - (j+1)\tilde{{\bf Y}}_{j+1}(\theta), \quad j=0, 1, \ldots, d.$$
This means that setting ${\bf y}_{k,d+1}={\bf 0}$ for $k \neq 0$, we can compute recursively  {for $j=d, d-1, \ldots, 1, 0$, }
$${\bf y}_{k,j} = \frac{1}{i k \omega} \left ({\bf z}_{k,j}- (j+1) {\bf y}_{k,j+1}\right ).$$

\item For $k=0$, one gets from \eqref{poly_q} and \eqref{p_and_hat_Z}
$$ {\bf y}_{0,j+1} = \frac{{\bf z}_{0,j}}{j}, \quad j=0, \ldots, d,$$
and,
$${\bf y}_{0,0} = {\bf y}_0 - \sum_{k=1}^M ({\bf y}_{k,0} + {\bf y}_{-k,0}).$$
\end{itemize}

\section{Taylor-Fourier methods and high order stroboscopic averaging}
\label{sec:averaging}

In this section, we establish a connection between semi-discrete Taylor-Fourier expansions and high-order stroboscopic averaging of systems of the form \eqref{eq:TF1_2_epsilon} based on the results in~\cite{Ch-M-SS2012a}. We begin by recalling from~\cite{Ch-M-SS2012a,Ch-M-SS2012} the series expansion \eqref{eq:y_epsilon} and its relation to high-order averaging, which, as shown in~\cite{Ch-M-SS2012a}, is derived through successive applications of the Picard-like iteration \eqref{eq:iter_1}--\eqref{eq:iter_3}. Note that, by rescaling time, we may assume without loss of generality that either $\omega=1$ or $\epsilon=1$ in \eqref{eq:TF1_2_epsilon}. However, we retain this redundancy in this section, as it facilitates the extension of our results to the quasi-periodic case considered in Subsection~\ref{subsec:quasi-periodic}.

In $d$th-order averaging, a change of variables of the form
\begin{equation*}
{\bf y} = {\bf U}(\omega t, {\bf w}; \epsilon) = {\bf w} + \epsilon \, {\bf u}_1(\omega t, {\bf w}) + \cdots + \epsilon^d \, {\bf u}_d(\omega t, {\bf w})
\end{equation*}
is introduced. This transformation converts the original system \eqref{eq:TF1_2_epsilon} into an averaged system of the form
\begin{equation*}
\frac{d}{dt}{\bf w} = \epsilon \, {\bf F}({\bf w}; \epsilon) + \epsilon^{d+1} {\bf R}(\omega t, {\bf w}; \epsilon), \quad {\bf w}(0) = {\bf w}_0,
\end{equation*}
where ${\bf F}({\bf w}; \epsilon)$ represents the averaged dynamics, and ${\bf R}(\omega t, {\bf w}; \epsilon)$ captures the remaining small amplitude oscillatory terms.

In the stroboscopic averaging framework, the solution ${\bf y}(t)$ of \eqref{eq:TF1_2_epsilon} admits the representation
\begin{equation}
\label{eq:averaging_approximation}
{\bf y}(t) =  {\bf U}(\omega t, {\bf W}(t,{\bf y}_0;\epsilon); \epsilon) + \mathcal{O}(\epsilon^{d+1}),
\end{equation}
where ${\bf W}(t,{\bf y}_0;\epsilon)={\bf w}(t)$ denotes the solution of the averaged system
\begin{equation}
\label{eq:averaged}
\frac{d}{dt}{\bf w} = \epsilon \, {\bf F}({\bf w}; \epsilon),  \quad {\bf w}(0)={\bf y}_0.
\end{equation}
Consequently, the solution ${\bf x}(t)$ of the original semi-linear equation
\begin{equation}
\label{eq:TF1_1_epsilon}
\frac{d}{dt}{\bf x} = \mathcal{A} {\bf x} + \epsilon\, {\bf g}({\bf x})
\end{equation}
satisfies
\begin{equation*}
{\bf x}(t) =  e^{t \mathcal{A}}\, {\bf U}(\omega t, {\bf W}(t,{\bf y}_0;\epsilon); \epsilon) + \mathcal{O}(\epsilon^{d+1}).
\end{equation*}

As shown in \cite[Theorem 2.10]{Ch-M-SS2012a} (see also \cite[Theorem 2.1]{Ch-M-SS2012}), we have that, for arbitrary ${\bf y}_0 \in \mathbb{R}^D$,
\begin{equation}
\label{eq:averaging_maps}
\begin{split}
{\bf U}(\omega t, {\bf y}_0;\epsilon) &=  {\bf y}_0 + \sum_{\ell=1}^{d} \epsilon^{\ell}
\sum_{k=-\infty}^{\infty} e^{i k \omega t} {\bf y}_{k,\ell,0}, \\
 {\bf W}(t,{\bf y}_0;\epsilon) &= {\bf y}_0 + \sum_{\ell=1}^{d} \epsilon^{\ell}   \sum_{j=0}^{\ell} t^j \sum_{k=-\infty}^{\infty} {\bf y}_{k,\ell,j}, \\
  {\bf F}({\bf y}_0; \epsilon) &= \frac{\partial}{\partial t}{\bf W}(0,{\bf y}_0;\epsilon) = \sum_{\ell=1}^{d} \epsilon^{\ell-1}  \sum_{k=-\infty}^{\infty} {\bf y}_{k,\ell,1},
 \end{split}
\end{equation}
where ${\bf y}_{k,\ell,j}$ are the coefficients of the expansion \eqref{eq:y_epsilon} of the solution of \eqref{eq:TF1_2_epsilon}. Recall that these coefficients can be obtained through successive applications of the Picard-like iteration \eqref{eq:iter_1}--\eqref{eq:iter_3}.

 The $(\infty,d)$-Taylor-Fourier approximation of the solutions ${\bf y}(t)$  and ${\bf x}(t)$
of \eqref{eq:TF1_2_epsilon} and  \eqref{eq:TF1_1_epsilon} read
\begin{equation}
\label{eq:Yd_epsilon}
{\bf Y}^{[d]}(\omega t, t; \epsilon, {\bf y}_0) =  \sum_{j=0}^{d} \epsilon^j t^j \, \sum_{k=-\infty}^{\infty} e^{i k \omega t} {\bf y}_{k,j}^{[d]}.
\end{equation}
and ${\bf X}^{[d]}(\omega t, t; \epsilon, {\bf y}_0) = e^{t\, \mathcal{A}}\, {\bf Y}^{[d]}(\omega t, t; \epsilon, {\bf y}_0)$, respectively.

Expansions \eqref{eq:Yd_epsilon} and \eqref{eq:y_epsilon_N} are derived through a similar Picard-like iteration process. The key distinction lies in the truncation method used to approximate ${\bf z}(t) = {\bf f}(\omega t, {\bf y}(t))$: \eqref{eq:Yd_epsilon} is obtained by truncating higher powers of $t$, whereas \eqref{eq:y_epsilon_N} involves truncating higher powers of $\epsilon$. This latter approach discards more terms in the expansion of ${\bf z}^{[d]}(t) = {\bf f}(\omega t, {\bf y}^{[d]}(t))$.  This implies that
the coefficients of the expansions \eqref{eq:Yd_epsilon} and \eqref{eq:y_epsilon_N} are related by
\begin{equation*}
\epsilon^j {\bf y}_{k,j}^{[d]} = \epsilon^j {\bf y}_{k,j,j} + \mathcal{O}(\epsilon^{d+1}).
\end{equation*}
We thus have the following relationships:
\begin{align*}
{\bf U}(\omega t, {\bf y}_0; \epsilon) &= {\bf Y}^{[d]}(\omega t, 0; \epsilon, {\bf y}_0) + \mathcal{O}(\epsilon^{d+1}), \\
{\bf W}(t, {\bf y}_0; \epsilon) &= {\bf Y}^{[d]}(0, t; \epsilon, {\bf y}_0)+ \mathcal{O}(\epsilon^{d+1}), \\
\epsilon\, {\bf F}({\bf y}_0; \epsilon) &= \left. \frac{\partial}{\partial t} {\bf Y}^{[d]}(0, t; \epsilon, {\bf y}_0) \right|_{t=0} + \mathcal{O}(\epsilon^{d+1}).
\end{align*}
This demonstrates that the change of variables in the $d$th-order stroboscopic averaging can be effectively approximated by applying the $(M,d)$-Taylor-Fourier algorithm with a sufficiently large $M$. In particular, this together with
\eqref{eq:averaging_approximation} implies that
\begin{equation}
\label{eq:averaging_approximation2}
{\bf X}^{[d]}\left(\omega t, 0; \epsilon, {\bf Y}^{[d]}(0, t; \epsilon, {\bf y}_0)\right)  =  {\bf X}^{[d]}(\omega t, t; \epsilon, {\bf y}_0)+ \mathcal{O}(\epsilon^d).
\end{equation}

While the right-hand side of the averaged equations and their solutions can also be computed using this approach, we do not claim it to be an efficient method for this purpose. However, it can serve as a useful tool to assess the accuracy of more efficient implementations of approximate averaged equations and their solutions.

We note that the semi-discrete expansions \eqref{eq:y_epsilon_N} in powers of $\epsilon$ can be transformed into a practical fully discrete method by adapting the procedure proposed in this paper. Specifically, for a given pair $(M,d)$, a finite-sum substitute for the $d$th order expansion \eqref{eq:y_epsilon_N} can be constructed by applying trigonometric interpolation to replace the periodic functions
\begin{equation*}
\sum_{k=-\infty}^{\infty} e^{i k \omega t} {\bf z}_{k,\ell,j}
\end{equation*}
in the first sub-step of the iteration defined by \eqref{eq:iter_1}--\eqref{eq:iter_3}.

This approach is particularly relevant in the context of stroboscopic averaging, especially when one aims to evaluate the mappings in \eqref{eq:averaging_maps} exactly --or with high precision for a sufficiently large $M$-- for a given order $d$. In particular, the right-hand side ${\bf F}({\bf y}_0; \epsilon)$ of the averaged equations obtained through this method preserves the geometric structure of the original equations. For instance, if the original equations \eqref{eq:TF1_1_epsilon} are Hamiltonian, then the averaged equations retain this Hamiltonian structure. In contrast, the approximation
\begin{equation*}
\left. \frac{\partial}{\partial t} {\bf Y}^{[M,d]}(0, t; \epsilon, {\bf y}0) \right|_{t=0} \approx {\bf F}({\bf y}_0; \epsilon),
\end{equation*}
obtained using the Taylor-Fourier algorithm proposed in Section~\ref{sec:implementation}, generally does not correspond to a Hamiltonian vector field. In this context, the truncated expansions \eqref{eq:y_epsilon_N} and \eqref{eq:Yd_epsilon} are not symplectic, as is expected for any explicit formula derived by truncating a power series expansion of a symplectic map.

However, this fully discrete method is computationally more expensive than the Taylor-Fourier algorithm proposed in this paper, both in terms of memory and complexity. For a given order $d$, it requires $(d+2)/2$ times more coefficients, and the algorithm involves arithmetic with truncated power series in two variables (rather than one).

\section{Application to more general equations}

\subsection{Semi-linear equations in the complex domain}

The definition of Taylor-Fourier approximations can be extended without modifications to initial value problems of the form  \eqref{eq:TF1_2} in the complex domain. In particular, Taylor-Fourier approximations can be successfully applied to equations of the form (\ref{eq:TF1_1}) defined in the complex domain, provided that  $A$ is a complex $D\times D$ matrix whose eigenvalues are integer multiples of the imaginary unit $i$, and the map
${\bf g}: \bC^D \to \bC^D$ is such that ${\bf g}({\bf x}) = {\bf G}({\bf x},{\bf \bar x})$, where ${\bf G}: \bC^{2D} \to \bC^D$ is a polynomial map.

\subsection{Semi-linear equations with non-resonant frequencies}
\label{subsec:quasi-periodic}

Taylor-Fourier approximations can be naturally extended to systems of the form \eqref{eq:TF0}, where $\mathcal{A}$ is a real $D \times D$ matrix with purely imaginary eigenvalues that are not fully resonant.
Applying the time-dependent change of variables
\begin{equation}
\label{eq:change_var_gen}
{\bf x}(t) = e^{t \mathcal{A}} {\bf y}(t),
\end{equation}
we transform \eqref{eq:TF0} into
\begin{equation*}
\frac{d}{dt}{\bf y} =   e^{-t \mathcal{A}} {\bf g}(e^{t \mathcal{A}} {\bf y}).
\end{equation*}

In this general setting, there exists a vector of non-resonant frequencies $\boldsymbol{\omega} = (\omega_1, \ldots, \omega_r) \in \mathbb{R}^r$ such that $\mathbf{k} \cdot \boldsymbol{\omega} = k_1 \omega_1 + \cdots + k_r \omega_r \neq 0$ for all $\mathbf{k} = (k_1, \ldots, k_r) \in \mathbb{Z}^r/\{{\bf 0}\}$.
The eigenvalues of $\mathcal{A}$ then take the form $(\mathbf{k} \cdot \boldsymbol{\omega}) i$ for $\mathbf{k} \in \mathbb{Z}^r$, making the solution operator $e^{t \mathcal{A}}$ of the linear part of \eqref{eq:TF0} quasi-periodic in $t$, rather than strictly periodic as in the fully resonant case.
This implies that there exists a mapping
\begin{equation*}
\begin{split}
{\bf f}: &\mathbb{T}^r \times \mathbb{R}^D \longrightarrow \mathbb{R}^D, \\
&\  (\boldsymbol{\theta}, {\bf y}) \mapsto {\bf f}(\boldsymbol{\theta}, {\bf y})
\end{split}
\end{equation*}
such that
\begin{equation*}
{\bf f}(\boldsymbol{\omega} t, {\bf y}) \equiv e^{-t \mathcal{A}} {\bf g}(e^{t \mathcal{A}} {\bf y}).
\end{equation*}
Therefore, the change of variables \eqref{eq:change_var_gen} transforms \eqref{eq:TF0} into the quasi-periodic system
\begin{equation}
\label{eq:TF0_y}
\frac{d}{dt}{\bf y} =  {\bf f}(\boldsymbol{\omega} t, {\bf y}).
\end{equation}

The definition of semi-discrete Taylor-Fourier approximations follows naturally: in Subsection~\ref{subsec:semi-discrete}, replace the basic frequency $\omega$ and the angular variable $\theta$ by the vector of non-resonant frequencies $\boldsymbol{\omega}$ and the angular variables $\boldsymbol{\theta}$, respectively.  The discussion on high order stroboscopic averaging in Section~\ref{sec:averaging} is also valid in the quasi-periodic case with the referred changes (the averaging technique is referred in that case as quasi-stroboscopic averaging~\cite{Ch-M-SS2012a,Ch-M-SS2012}).

Fully discrete $(M,d)$-Taylor-Fourier approximations extend this approach by employing trigonometric interpolation of functions $\mathbb{T}^r \to \mathbb{C}$ over an $r$-fold Cartesian product of the uniform grid
\begin{equation*}
\theta_n = \frac{n\pi}{M}, \quad n=0, 1, 2, \ldots, 2M-1,
\end{equation*}
in $[0,2\pi)$, which can be efficiently computed using multi-dimensional discrete Fourier transforms (DFT).

However, compared to the periodic case, quasi-periodicity requires increased memory and higher computational cost: while a complex tensor of dimensions $D \times (d+1) \times (M+1)$ suffices to store the coefficients ${\bf y}^{[M,d]}_{k,j}$ in the periodic case, the quasi-periodic setting with $\boldsymbol{\omega} \in \mathbb{R}^r$ requires a complex tensor of dimensions $D \times (d+1) \times (M+1)^r$, significantly increasing storage and processing demands.

\subsection{Equations in other alternative formats}
Going back to equations in the real domain,  Taylor-Fourier approximations are not limited to semi-linear equations of the form (\ref{eq:TF1_1}).
They can be applied more generally to problems of the form
\begin{equation}
\label{eq:TF_gen_ode}
\frac{d}{dt}{\bf x} = \omega\, {\bf r}( {\bf x}) + {\bf g}({\bf x}), \quad {\bf x}(0)= {\bf x}_0,
\end{equation}
where ${\bf r}: \bR^D \to \bR^D$ and   ${\bf g}: \bR^D \to \bR^D$ are smooth maps and $\omega$ is a positive real parameter,
provided that the following two assumptions hold:
\begin{enumerate}
\item  {For each $t\in\bR$ there exists} a $t$-flow map $\varphi_t:\mathbb{R}^D \to \mathbb{R}^D$ of the differential system $\frac{d}{dt}{\bf x} = {\bf r}( {\bf x})$. That is, for each $t \in \bR$ and ${\bf x} \in \bR^D$,
\begin{equation}
\frac{\partial}{\partial t} \varphi_t({\bf x}) = {\bf r}( \varphi_t({\bf x})).
\end{equation}
 {In that case, problem} (\ref{eq:TF_gen_ode}) is transformed with the change of variables
\begin{equation*}
{\bf x}(t) = \varphi_{\omega t}( {\bf y}(t))
\end{equation*}
into (\ref{eq:TF1_2}), where
\begin{equation}
\label{eq:f_gen}
{\bf f}(\theta, {\bf y}) = \varphi'_\theta({\bf y})^{-1} {\bf g}(\varphi_\theta({\bf y})).
\end{equation}
Here, $\varphi'_{\theta}({\bf y})$ denotes the Jacobian matrix of the map $\varphi_{\theta}$ evaluated at the vector ${\bf y}$.
\item The  {function ${\bf f}(\theta, {\bf y})$ defined in} (\ref{eq:f_gen}) is $(2\pi)$-periodic in $\theta$.
\end{enumerate}
 {Under these two assumptions}, the solution ${\bf x}(t)$ of \eqref{eq:TF_gen_ode} can be approximated by
\begin{equation}
\label{eq:xd_2}
{\bf x}^{[M,d]}(t) = \varphi_{\omega \, t}\left({\bf y}^{[M,d]}(t) \right),
\end{equation}
where
\begin{equation*}
{\bf y}^{[M,d]}(t) = \sum_{k=-M}^M e^{ik\omega t} \sum_{j=0}^{d} t^j \, {\bf y}_{k,j}^{[M,d]}
\end{equation*}
is the $(M,d)$-Taylor-Fourier approximation  {to} the solution of \eqref{eq:TF1_2} with the map ${\bf f}$ defined by \eqref{eq:f_gen}.

In particular, the above-mentioned two assumptions hold for systems of the form
\begin{equation}
\label{eq:aa}
\begin{split}
\frac{d}{dt} \phi_j &= k_j \omega + g_j({\bf x}), \quad j=1,\ldots,s, \\
\frac{d}{dt} a_j &=  g_{s+j}({\bf x}), \quad j=1,\ldots,D-s, \\
\end{split}
\end{equation}
where ${\bf x} = (\phi_1,\ldots,\phi_{s},a_1,\ldots,a_{D-s})$,  $k_j \in \mathbb{Z}$,  and for each $j$, $g_j$ is  $2\pi$-periodic in each of the angle variables $\phi_j$.  In that case,
\begin{equation*}
\varphi_{t}(\phi_1,\ldots,\phi_{s},a_1,\ldots,a_{D-s}) =
\left(
\begin{matrix}
\phi_1 + k_1 \, t\\
\vdots\\
\phi_{s} + k_{s}\, t\\
a_1\\
\vdots \\
a_{D-s}\\
\end{matrix}
\right),
\end{equation*}
and the change of variables ${\bf x}(t) = \varphi_{\omega t}( {\bf y}(t))$, with  $ {\bf y} = (\tilde \phi_1,\ldots,\tilde \phi_{s},\tilde a_1,\ldots,\tilde a_{D-s})$,
transforms the system (\ref{eq:aa}) into
\begin{align*}
\frac{d}{dt} \tilde \phi_j &= g_j(\tilde \phi_1+ k_1 \, \omega \,t,\ldots,\tilde \phi_{s}+ k_s\, \omega \,t,\tilde a_1,\ldots,\tilde a_{D-s}), \quad j=1,\ldots,s, \\
\frac{d}{dt} \tilde a_j &=  g_{s+j}(\tilde \phi_1+ k_1 \, \omega \,t,\ldots,\tilde \phi_{s}+ k_s \, \omega \,t,\tilde a_1,\ldots,\tilde a_{D-s}), \quad j=1,\ldots,D-s.
\end{align*}

\section{Numerical illustrations}
\label{sec:experiments}
In this section we include two numerical examples to illustrate the performance of the Taylor-Fourier approximations proposed
in this paper.  All codes we have used to carry out the numerical experiments are available in  \url{https://github.com/jmakazaga/Taylor-Fourier}.

\subsection{The cubic nonlinear Schr\"odinger equation}
\label{sec:NLS}
The first example is the cubic nonlinear Schr\"odinger (NLS) equation
\begin{equation}
\label{eq:NLS}
i u_t + u_{xx} + |u|^2 u = 0, \quad x \in [0, 2\pi], 
\end{equation}
with periodic boundary conditions and prescribed initial data
\begin{equation}
\label{NLS_gic}
u(x,0)=u^0(x).
\end{equation}
For the numerical solution of \eqref{eq:NLS}-\eqref{NLS_gic}, we use the method of lines and then, we first discretize in space with the spectral collocation method as follows.

We fix an integer $J$, set $h=\pi/J$, introduce the uniform grid of points $x_j = (j-1)h, 1 \leq j \leq 2J$, and denote
${\bf U}(t) = [U_1(t), \ldots, U_{2J}(t)]^T$, where $U_j(t)$ is an approximation to $u(x_j,t)$, for $1 \leq j \leq 2J$. We also consider the vector $\hat{\bf U}(t) = [\hat{U}_{1}(t), \ldots, \hat{U}_{2J}(t)]^T$ obtained by applying the discrete Fourier transform to ${\bf U}(t)$, that is, for $k=1, \ldots, 2J$,
$$
\hat{U}_k(t) = \frac{1}{2J} \sum_{j=1}^{2J} U_j(t) e^{-i(k-1)x_j}.
$$
In matrix notation, we write it as $\hat{\bf U} = \calF {\bf U}$. The vector ${\bf U}$ can be recovered from
$\hat{\bf U}$ by applying the inverse Fourier transform, i.e., ${\bf U} = \calF^{-1} \hat{\bf U}$. More explicitly, for $j=1, \ldots, 2J$,
$$U_j(t) =  \frac{1}{2J} \sum_{k=1}^{2J} \hat{U}_k(t) e^{i(k-1)x_j},$$
or equivalently,
$$U_j(t) =  \sum_{k=-J}^{J} \hat{u}_k(t) e^{ikx_j},$$
where
\begin{equation}
\begin{array}{l}
\hat{u}_0(t) = \dfrac{1}{2J} \hat{U}_1(t), \quad \hat{u}_J(t) = \hat{u}_{-J}(t) = \dfrac{1}{4J} \hat{U}_{J+1}(t), \\[8pt]
\hat{u}_k(t) = \dfrac{1}{2J} \hat{U}_{k+1}(t), \quad \hat{u}_{-k}(t) = \dfrac{1}{2J} \hat{U}_{2J-k+1}(t), \quad k=1, \ldots, J-1.
\end{array}
\end{equation}
Notice that
$$\tilde{u}(x,t) = \sum_{k=-J}^J \hat{u}_k(t)e^{ikx}$$
is a trigonometric polynomial satisfying the interpolation conditions
$$\tilde{u}(x_j,t) = U_j(t), \quad j=1, \ldots, 2J.$$

Both, the discrete Fourier transform and its inverse, can be efficiently computed by applying the fast Fourier transform (FFT) algorithm. This is particularly efficient when $J = 2^{l}$.

The spectral collocation semi-discretization of the nonlinear Schr\"odinger equation is obtained by requiring that
$$
i \tilde{u}_t(x_j, t) + \tilde{u}_{xx}(x_j,t) + |\tilde{u}(x_j,t)|^2 \tilde{u}(x_j,t) = 0, \quad  \quad j=1, \ldots, 2J.
$$
By taking into account that
\begin{eqnarray*}
i \, \tilde{u}_{xx}(x_j,t) &=& \sum_{k=-J}^J -i \, k^2 \, \hat{u}_k(t) e^{ikx_j} \\
&=& \frac{1}{2J}\sum_{k=1}^J \lambda_k \, \hat{U}_k(t) e^{i(k-1)x_j},
\end{eqnarray*}
where
$$\lambda_k = \left \{ \begin{array}{rcl} -i(k-1)^2 &\mbox{ for }& k=1, \ldots, J+1, \\ -i(2J-k+1)^2 &\mbox{ for }& k=J+2, \ldots, 2J,\end{array} \right. $$
we arrive at the semi-discretized NLS equation
\begin{equation}
\label{sd_NLS}
\frac{d}{dt} {\bf U} = A {\bf U} + {\bf g}({\bf U}),
\end{equation}
where
\begin{equation}
\label{g_NLS}
{\bf g}({\bf U}) = i \, \left ( \begin{matrix}|U_1|^2 U_1 \\ \vdots \\ |U_{2J}|^2 U_{2J}\end{matrix} \right ),
\end{equation}
$A = {\calF}^{-1} D {\calF}$, and $D$ is the diagonal matrix with diagonal entries
$$(\lambda_1, \ldots, \lambda_{2J}) = (0, -i, -2^2 i, \ldots, -(J-1)^2 i, -J^2i, -(J-1)^2 i, \ldots, -2^2 i, -i).$$

The time-dependent change of variables ${\bf U}=e^{tA}{\bf W}$ transforms (\ref{sd_NLS}) into
\begin{equation}
\label{tsd-NLS}
\frac{d}{dt} {\bf W} = {\bf f}(t, {\bf W}),
\end{equation}
where ${\bf f}(t, {\bf W})$ is $(2\pi)$-periodic in $t$ and is defined by
$${\bf f}(\theta, {\bf W}) = e^{-\theta A} {\bf g}\left (e^{\theta A}{\bf W} \right ).$$
Products of a vector in $\bC^{2J}$ by matrices $e^{\theta A}$ and $e^{-\theta A}$ can be computed efficiently by performing two FFT in $\bC^{2J}$, since $e^{\theta A} = {\cal F}^{-1}e^{\theta D} {\cal F}$, and $e^{\theta D}$ is the diagonal matrix with diagonal entries $(e^{\lambda_1 \theta}, \ldots, e^{\lambda_{2J} \theta})$. On the other hand, due to the particular form of the nonlinearity ${\bf g}$ in \eqref{g_NLS}, the Taylor expansions required to compute the Taylor-Fourier approximations are cheaply available since sums and products of power series in $t$ can be computed in a very efficient way.

In the numerical experiments we have considered the NLS equation with initial data $u(x,0)=\eta(x)$ equal to the step function
\begin{equation}
\label{NLS_ic}
\eta(x)=\left \{ \begin{array}{l} -1 \quad \mbox{ if } x \in [0, \pi), \\
\hphantom{-}1 \quad \mbox{ if } x \in [\pi, 2\pi], \end{array} \right.
\end{equation}
a rough initial data that has widely appeared in the recent literature~\cite{Ostermann-Schratz,Ostermann-Yao, Chen-Olver}. For instance, in \cite{Chen-Olver} the initial condition \eqref{NLS_ic} is used to illustrate the Talbot effect.

\begin{figure}[H]
\centering
{
\includegraphics[width=0.45\hsize]{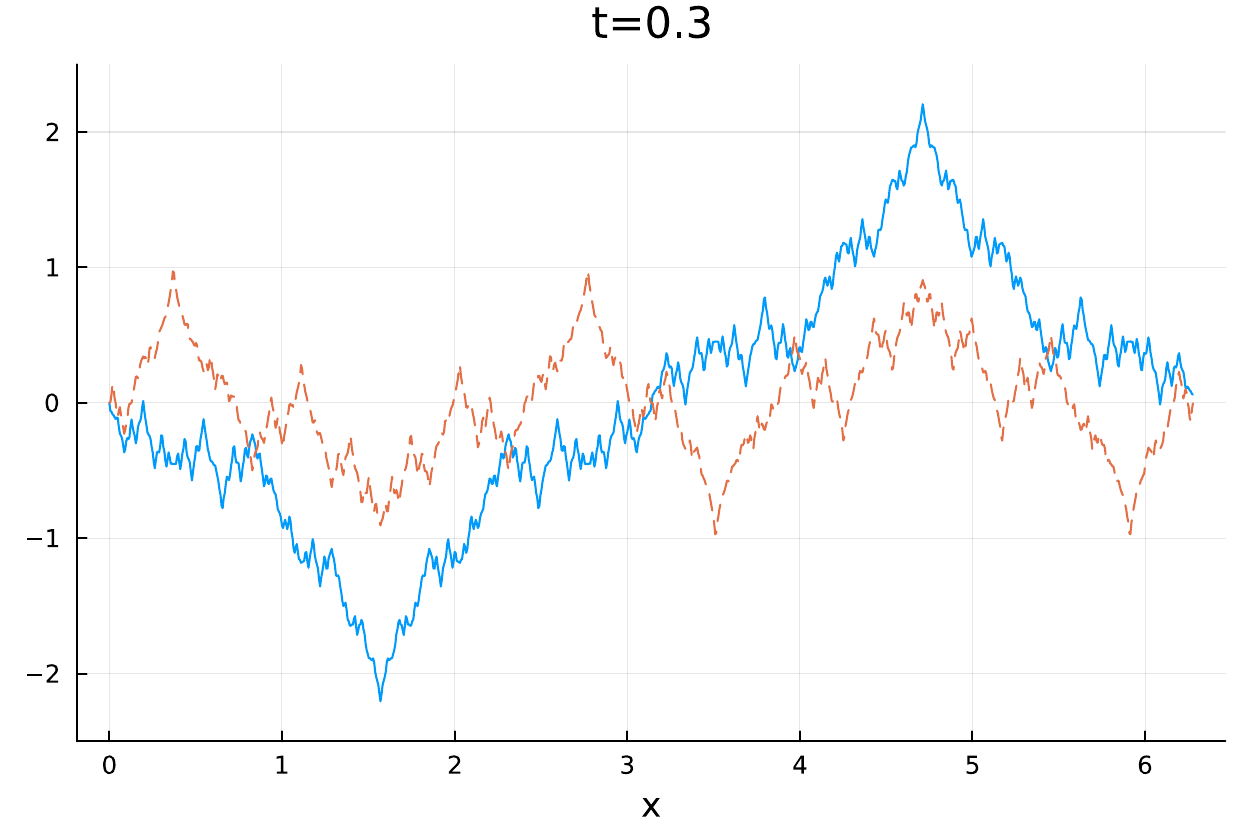}\quad
\includegraphics[width=0.45\hsize]{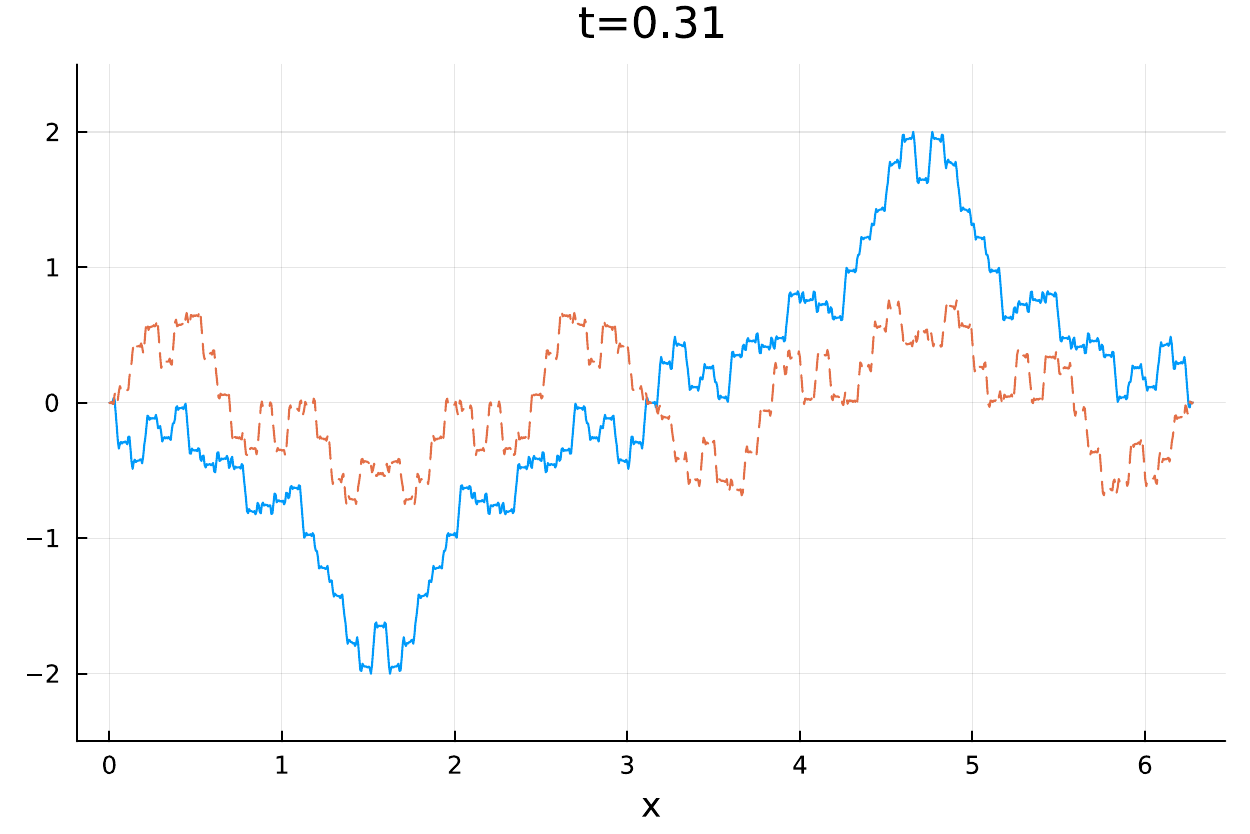} \\
\includegraphics[width=0.45\hsize]{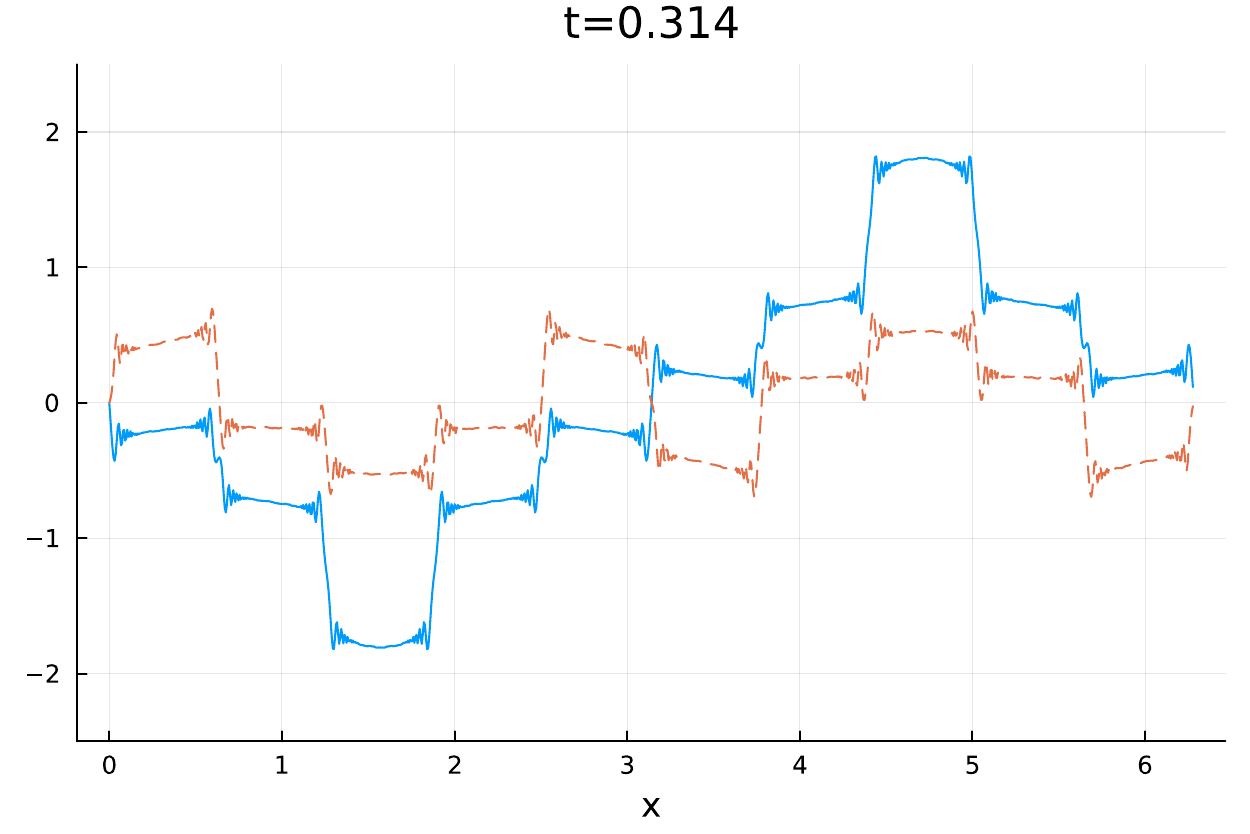}\quad
\includegraphics[width=0.45\hsize]{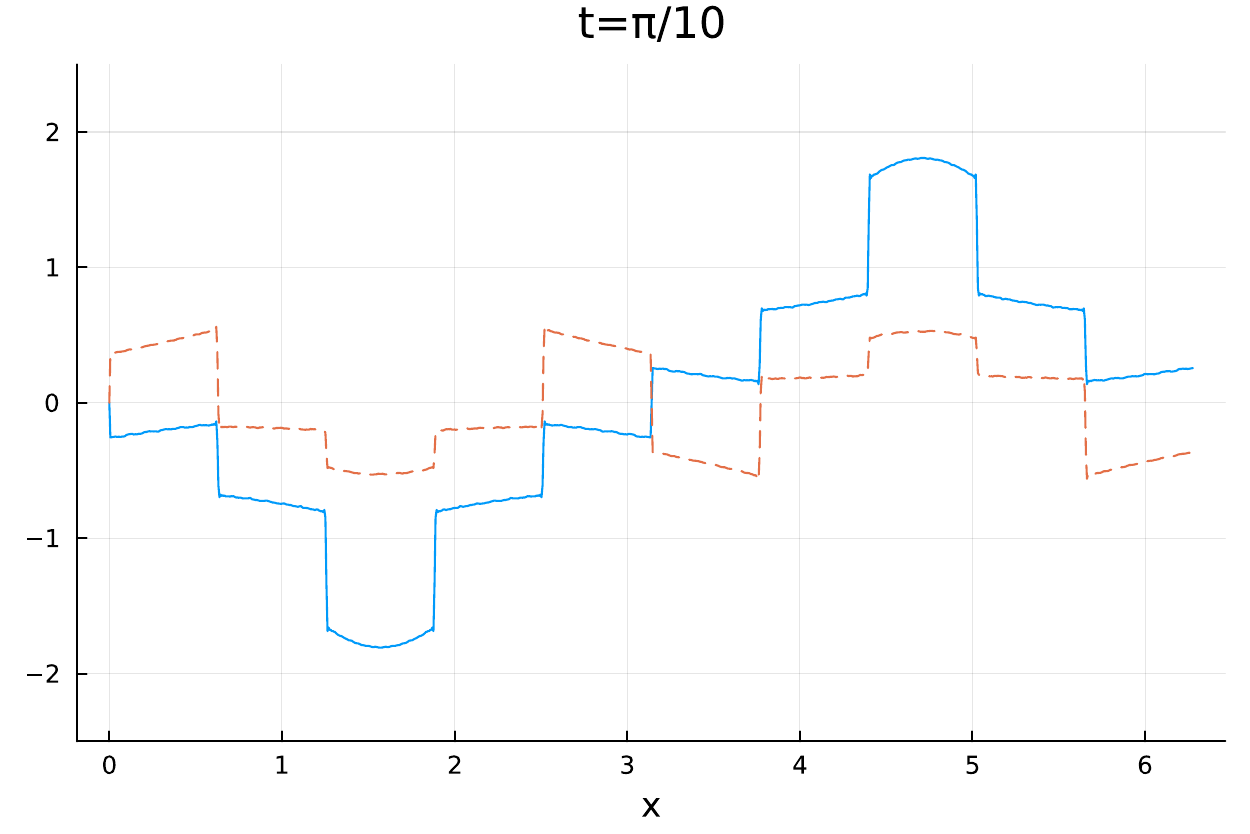}
}
\caption{Numerical solution of \eqref{eq:NLS}-\eqref{NLS_ic} at $t=0.3$ (top left), $t=0.31$ (top right), $t=0.314$ (bottom left) and $t=\pi/10$ (bottom right) computed with $J=2^9$,  {and the $(2^{11},3)$-Taylor-Fourier approximation}. Real part is plotted with blue solid line and imaginary part with brown dashed line.}
\label{fig:NLS1}
\end{figure}

In the first  experiment, we take $J=2^7$, $2^8$, $2^9$, $2^{10}$, $T=\pi/10$, $d \le 5$, and for each $J$, $M=J$, $2J$, $4J$. Once computed the coefficients of the  $(M,d)$-Taylor-Fourier approximation, we have evaluated ${\bf W}^{[d]}(t)$ at $t=0.3$, $0.31$, $0.314$ and $\pi/10$, as in \cite{Chen-Olver}. The results for $J=2^9$,  $M=4J$ and $d=3$ are shown in Figure \ref{fig:NLS1}, to be compared with those in Figures 5 and 6 in \cite{Chen-Olver}. Different from \cite{Chen-Olver}, we have included in the same plot the real and the imaginary part of the numerical solutions.

\begin{figure}[H]
\centering
{
\includegraphics[width=0.30\hsize]{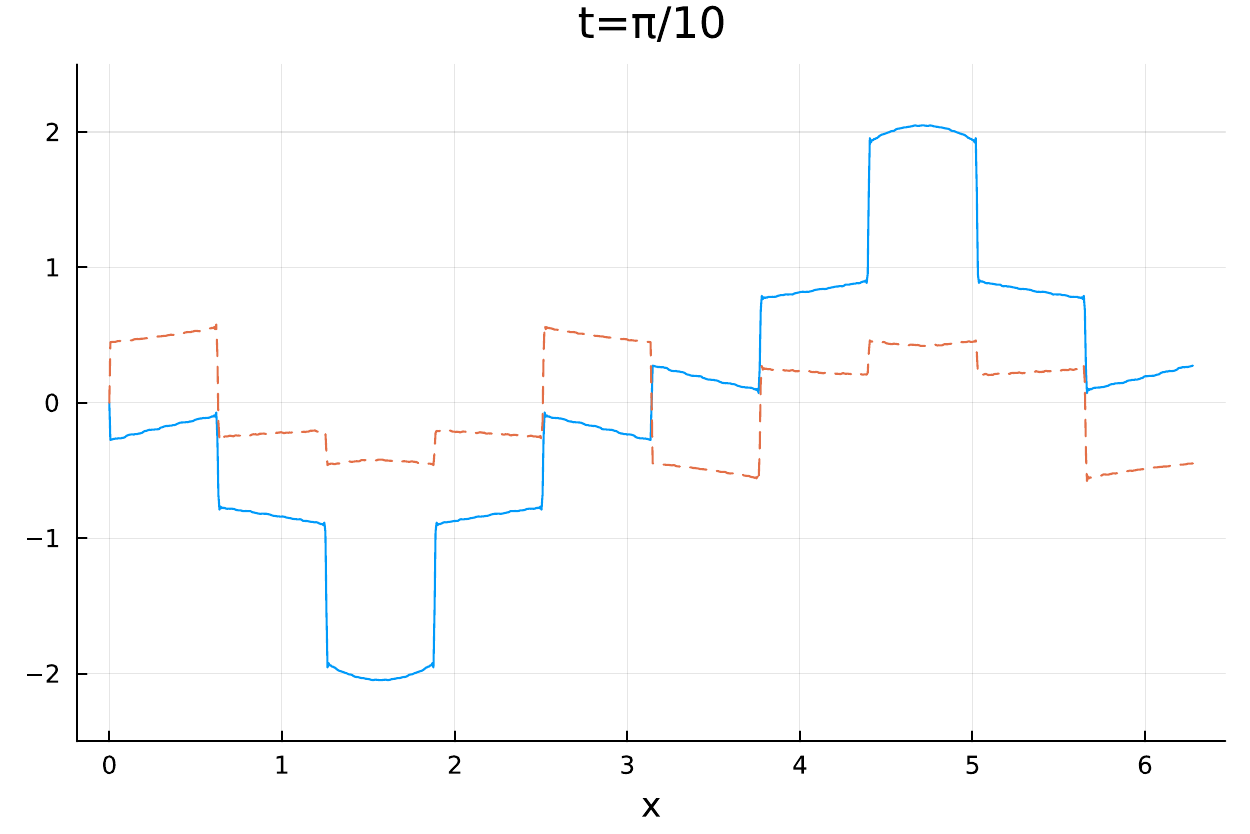} \quad
\includegraphics[width=0.30\hsize]{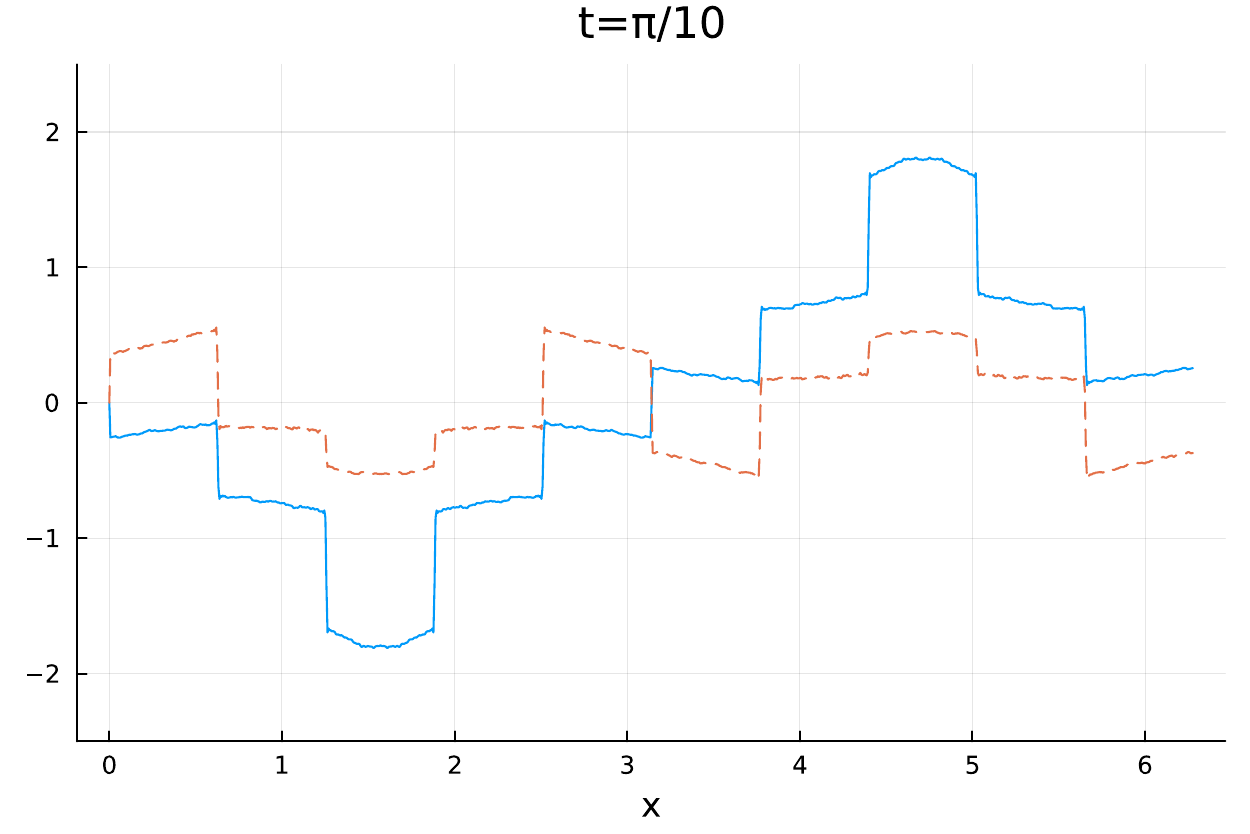} \quad
\includegraphics[width=0.30\hsize]{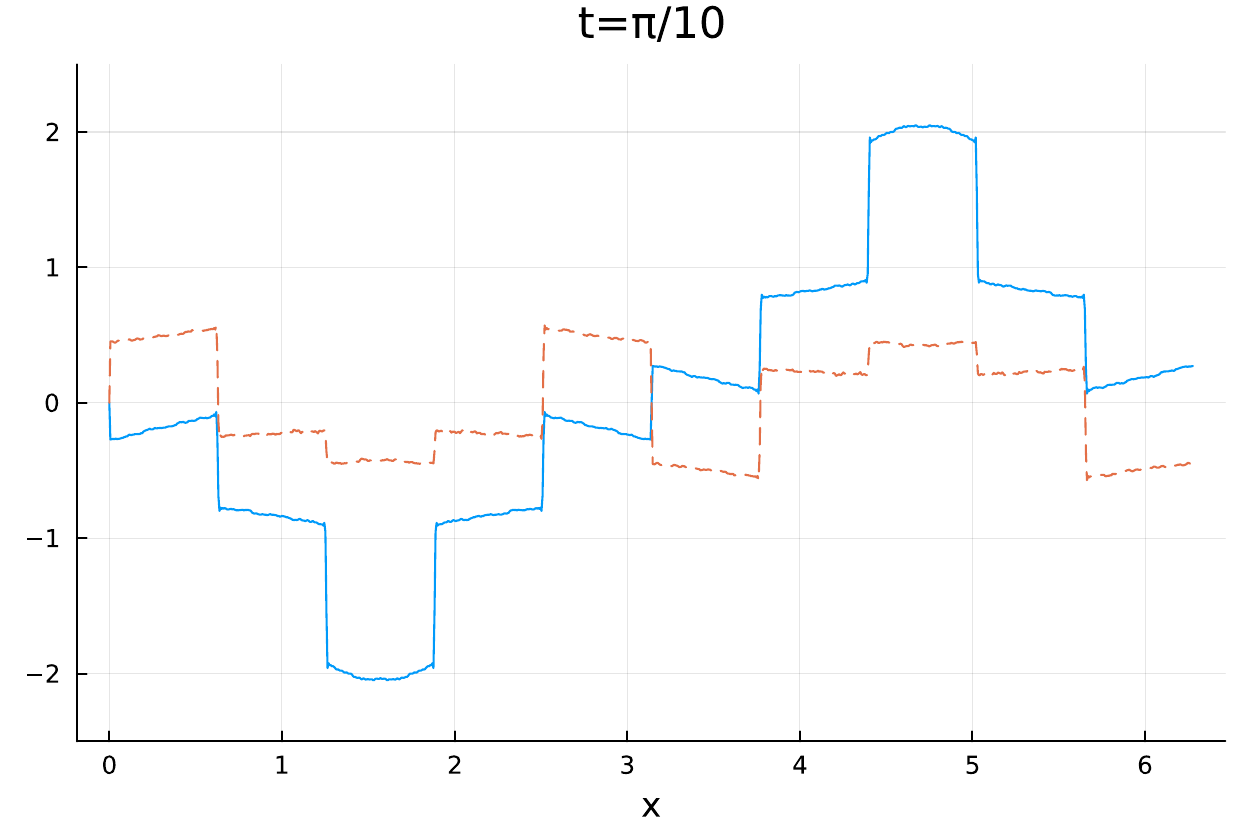}
}
\caption{ $(M,d)$-Taylor-Fourier approximations to the solution of \eqref{eq:NLS}-\eqref{NLS_ic} at $t=\pi/10$ for $J=2^9$,  computed with  $M=4J$, $d=1$ (left), $M=J$, $d=3$ (center) and $M=J$, $d=1$ (right). Real part is plotted with blue solid line and imaginary part with brown dashed line.}
\label{fig:NLS2}
\end{figure}

Secondly, to illustrate how the choice of $d$ and $M$ affects the approximate solution computed with the Taylor-Fourier method, we include in Figure~\ref{fig:NLS2} the numerical approximations to the solution of \eqref{eq:NLS}-\eqref{NLS_ic} at $t=\pi/10$ when $J=2^9$, and the combinations
 $M=4J$, $d=1$ (left), $M=J$, $d=3$ (center) and $M=J$, $d=1$ (right) are used. Plots in Figure~\ref{fig:NLS2} must be compared with the bottom right plot in Figure~\ref{fig:NLS1} where  {the $(4J,3)$-Taylor-Fourier approximation} is displayed. Comparing the plots corresponding to the same value of $M$ and reducing $d$ from 3 to 1, we observe that, although the overall shape of the numerical solutions is rather similar, the maximum and minimum values of both solutions are clearly different. Looking now at the solutions obtained with the same value of $d$ we see that when decreasing $M$ from $4J$ to $J$ small oscillations appear.

In order to test the effect of changing the basic frequency $\omega$, we consider (\ref{eq:NLS}) with initial data
\begin{equation}
\label{eq:initial_data_NLS}
u(0,x) = \epsilon\, \eta(x)
\end{equation}
and $t \in [0,\epsilon^{-2} \pi/10]$, with different values of $\epsilon>0$. This is equivalent to considering $u(t,x)=\epsilon\, v(\epsilon^2\, t, x)$, where $v(\tau,x)$ is the solution of
\begin{equation*}
i \, v_{\tau} = -\omega\, v_{xx} - |v|^2 v, \quad x \in [0, 2\pi], \quad \tau \in [0, \pi/10],
\end{equation*}
with $\omega=\epsilon^{-2}$.

\begin{figure}[H]
\centering
\includegraphics[width=0.8\hsize]{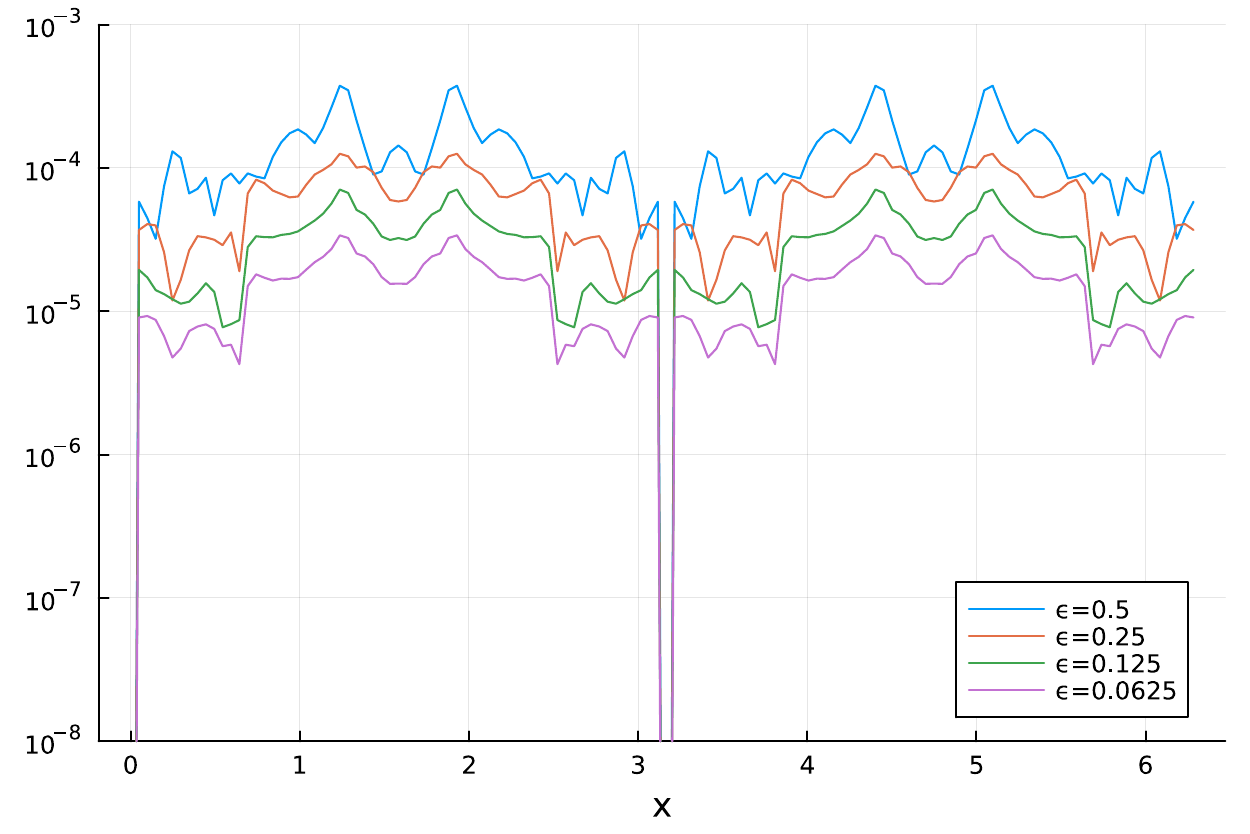}
\caption{Errors of the $(2^{11},7)$-Taylor-Fourier approximation at $t=\epsilon^{-2} \pi/10$ for $\epsilon=2^{-m}$, $m=1,2,3,4$ with $J=2^6$ and initial data (\ref{eq:initial_data_NLS}).}
\label{fig:NLS4}
\end{figure}

In the abstract and in Section~1 we claim that the accuracy of Taylor-Fourier approximations does not deteriorate as the basic frequency $\omega$ increases. We next check whether this holds true in the particular case of the ODE system obtained after the space discretization with $J=2^6$ of (\ref{eq:NLS}) for the $(2^{11},7)$-Taylor-Fourier approximation, with initial data (\ref{eq:initial_data_NLS}) for different values of $\epsilon>0$. In Figure~\ref{fig:NLS4}, we represent (in logarithmic scale) the error $|\tilde U_j(t)-U_j(t)|$ ($j=1,\ldots,2J$) at time $t=\epsilon^{-2} \pi/10$ versus $x_j$, $j=1,\ldots,2J$ for  $\epsilon=2^{-m}$, $m=1,2,3,4$.
{Here, $\tilde U_j(t)$ denotes the reference solution of the semi-discretized system
computed using the Lawson integrator based on the implementation DP5  in the julia library DifferentialEquation.jl of the 5th order embedded explicit Runge-Kutta pair proposed  in~\cite{DP5}, with tolerance for absolute and relative errors equal to $10^{-12}$.
 We observe that the errors do not increase as $\omega=\epsilon^{-2}$ increases, but rather decrease linearly with $\epsilon = \omega^{-1/2}$ in this case.}

\subsection{A perturbed Kepler problem}
\label{sec:satellite}
The second example is a perturbed Kepler problem (the main problem of satellite theory), which accounts for the oblateness of the earth. More precisely, it is the Kepler problem perturbed by the $J_2$-term,
\begin{align}
\label{eq:pkepler}
\ddot {\bf q} &= -\frac{\mu}{r^3} {\bf q} - \nabla V({\bf q}), \\
\label{eq:icond_pkepler}
 {\bf q}(t_0) &= {\bf q}_0, \quad \dot {\bf q}(t_0) = \dot {\bf q}_0,
\end{align}
where
\begin{equation*}
{\bf q} = \left(
\begin{matrix}
x\\
y\\
z
\end{matrix}
\right), \quad
r = \|{\bf q}\| = \sqrt{x^2+y^2+z^2}, \quad \dot {\bf q} = \frac{d}{dt} {\bf q}, \quad \ddot {\bf q} = \frac{d^2}{dt^2} {\bf q},
\end{equation*}
and
\begin{equation}
\label{eq:J2potential}
V({\bf q}) = J_2\, \frac{\mu\, R_e^2}{2\, r^3} \,\left(3\, \sin^2\theta -1 \right), \quad \sin\theta = \frac{z}{r}.
\end{equation}
Here, $J_2$ is the zonal harmonic coefficient ($J_2=1.08262668 \times 10^{-3}$), $\mu$ is the standard gravitational parameter of the Earth ($\mu = 398600.44189 \, km^3/seg^{2}$) and $R_e$ is the equatorial radius of the Earth ($R_e = 6378.137\, km$).

From now on, we will use the notation $\epsilon:= J_2\, \mu\, R_e^2$. From (\ref{eq:J2potential}), we get that
\begin{equation*}
\nabla V(x,y,z) = \frac{3\, \epsilon}{2\, r^5}
\left(
\begin{matrix}
(1 - 5\, \sin^2\theta)\, x\\
(1 - 5\, \sin^2\theta)\, y \\
(3 - 5\, \sin^2\theta)\, z
\end{matrix}
\right).
\end{equation*}


\subsubsection{KS-formulation of the equations of motion}
\label{s:KS}

We will rewrite the equations by making use of the Kustaanheimo-Stiefel (KS) transformation.

The position variables ${\bf q} = (x,y,z) \in \R^3$ are replaced by new variables ${\bf u}=(u_1, u_2, u_3,u_4) \in \R^4$ satisfying the relation ${\bf q} = L({\bf u}) {\bf u}$,
where
\begin{equation*}
L({\bf u}) =
\left(
\begin{matrix}
u_1 & -u_2 & -u_3 & u_4 \\
u_2 & u_1 & -u_4 & -u_3 \\
u_3 & u_4 & u_1 & u_2 \\
\end{matrix}
\right).
\end{equation*}
Clearly, it holds that $\|{\bf u}\|^2 = \|{\bf q}\| = r$.

Instead of the physical time $t$, a new independent variable (the {\em fictitious time})
\begin{equation*}
\tau = \int_0^t \frac{1}{\|{\bf q}(t)\|} \, dt
\end{equation*}
will be used (we will denote with a prime the derivative with respect to $\tau$), and the physical time $t$ is considered as an additional dependent variable of the equations of motion.

Following~\cite{Stiefel-Scheifele}, the solution curve ${\bf q}(t)$ of the initial value problem (\ref{eq:pkepler})--(\ref{eq:icond_pkepler}) can be parametrized as ${\bf q} = L({\bf u}(\tau)) {\bf u}(\tau)$, $t = t(\tau)$, where $({\bf u}(\tau),t(\tau))$ is the solution of the initial value problem
\begin{align}
\label{eq:odeu}
{\bf u}'' &= -\frac{h}{2} {\bf u} -  \nabla R({\bf u}), \\
\label{eq:u0}
{\bf u}(0)&={\bf u}_0, \quad {\bf u}'(0) = {\bf u}'_0, \\
\label{eq:odet}
t' &= \|{\bf u}\|^2, \\
\label{eq:t0}
t(0)&=t_0,
\end{align}
where
$$
h = \frac{\mu}{\|{\bf q}_0\|} - \frac12\, \langle \dot {\bf q}_0, \dot {\bf q}_0 \rangle - V({\bf q}_0),$$ and
\begin{equation*}
R({\bf u}) = \frac12\,  \|{\bf u}\|^2\, V(L({\bf u}) {\bf u}) = \frac{\epsilon}{4\, \|{\bf u}\|^4} \,\left(3\, \sin^2\theta -1 \right),
\end{equation*}
with
$\sin\theta = 2\, (u_1 u_3 + u_2 u_4)/\|{\bf u}\|^2$.

The initial condition ${\bf u}_0 \in \R^4$ is chosen in such a way that
\begin{equation}
\label{eq:q0}
{\bf q}_0 = L({\bf u}_0) {\bf u}_0,
\end{equation}
and ${\bf u}_0'\in \R^4$ is determined as
\begin{equation}
\label{eq:uprime0}
{\bf u}_0' = \frac12 \, L({\bf u}_0)^T \dot {\bf q}_0.
\end{equation}
There are infinitely many choices of ${\bf u}_0 \in \R^4$ satisfying (\ref{eq:q0}) for a given ${\bf q}_0 = (x_0,y_0,z_0)$. Among them, we will uniquely determine ${\bf u}_0 
= (u_{1,0}, u_{2,0}, u_{3,0}, u_{4,0})$ by imposing the following conditions:
\begin{equation}
\label{eq:xi}
\begin{array}{ll}
\mbox{if } x_0\geq 0, \quad u_{1,0} = u_{4,0} \geq 0, \\
\mbox{if } x_0<0, \quad u_{2,0} = u_{3,0} \geq 0.
\end{array}
\end{equation}
This leads to
$$u_{1,0} = u_{4,0} = \frac12\, \sqrt{r_0 + x_0}, \quad
u_{2,0} = \frac{y_0 u_{1,0} + z_0 u_{4,0}}{r_0 + x_0}, \quad u_{3,0} = \frac{z_0 u_{1,0} - y_0 u_{4,0}}{r_0 + x_0},$$
if $x_0\geq 0$, and
$$u_{2,0} = u_{3,0} = \frac12\, \sqrt{r_0 - x_0}, \quad
u_{1,0} = \frac{y_0 u_{2,0} + z_0 u_{3,0}}{r_0 - x_0}, \quad u_{4,0} = \frac{z_0 u_{2,0} - y_0 u_{3,0}}{r_0 - x_0},$$
if $x_0< 0$, where $r_0 = \|{\bf q}_0\|$.

The gradient of $R({\bf u})$ can be computed as
\begin{equation*}
\nabla R({\bf u}) =  \
\frac{\epsilon}{ \|{\bf u}\|^6} \, (1-6 \sin^2 \theta)
\left(
\begin{matrix}
u_1\\
u_2\\
u_3\\
u_4
\end{matrix}
\right)
+\frac{3\, \epsilon}{ \|{\bf u}\|^6} \, \sin \theta
\left(
\begin{matrix}
u_3\\
u_4\\
u_1\\
u_2
\end{matrix}
\right).
\end{equation*}
Clearly, the system of second order ODEs (\ref{eq:odeu}) is of the form (\ref{eq:TF1_1}) with
\begin{equation*}
\omega = \sqrt{h/2}, \quad
{\bf x} =
\left(
\begin{matrix}
{\bf u}\\
{\bf u}'
\end{matrix}
\right), \quad
A =
\left(
\begin{matrix}
0 & \omega^{-1} I_4\\
-\omega\, I_4  & 0
\end{matrix}
\right),
\quad
g({\bf x}) =
\left(
\begin{matrix}
{\bf 0}\\
-\nabla R({\bf u})
\end{matrix}
\right),
\end{equation*}
where $I_4$ denotes the $4\times4$ identity matrix.

\subsubsection{Stiefel and Scheifele's formulation}

From now on, we assume that $h>0$, and denote $\omega=\sqrt{h/2}$.  The solution ${\bf u}(\tau)$ of the initial value problem (\ref{eq:odeu})--(\ref{eq:u0}) can be written as~\cite{Stiefel-Scheifele}
\begin{equation*}
{\bf u}(\tau) = \cos(\omega \tau) \balpha(\tau)+ \omega^{-1} \sin(\omega \tau) \bbeta(\tau),
\end{equation*}
where $\balpha(\tau)$ and $\bbeta(\tau)$ are the solutions of the 8-dimensional ODE system
\begin{equation}
\label{eq:ode_alphabeta}
\begin{split}
\balpha' &= \omega^{-1} \sin(\omega \tau) \nabla R( \cos(\omega \tau) \balpha+ \omega^{-1} \sin(\omega \tau) \bbeta), \\
\bbeta' &= -\cos(\omega \tau)  \nabla R( \cos(\omega \tau) \balpha+ \omega^{-1} \sin(\omega \tau) \bbeta),
\end{split}
\end{equation}
supplemented with the initial conditions
\begin{equation}
\label{eq:alphabeta_icond}
\balpha(0) = {\bf u}_0,  \quad
\bbeta(0) = {\bf u}'_0.
\end{equation}
Actually,  (\ref{eq:ode_alphabeta}) is the VOP formulation (\ref{eq:TF1_2}) corresponding to  (\ref{eq:odeu}).

The solution $t(\tau)$ of (\ref{eq:odet})--(\ref{eq:t0}) can be  {written} as
\begin{equation*}
t(\tau) = t_0 + \int_0^{\tau} \| \cos(\omega \sigma) \balpha(\sigma)+ \omega^{-1} \sin(\omega \sigma) \bbeta(\sigma)\|^2\, d\sigma,
\end{equation*}
which can be explicitly computed once appropriate  $(M,d)$-Taylor-Fourier approximations
\begin{equation*}
\balpha(\tau) \approx \sum_{k=-M}^M  \sum_{j=0}^{d} e^{ik\omega \tau}\, \tau^j \, \balpha_{k,j},
\quad
\bbeta(\tau) \approx \sum_{k=-M}^M  \sum_{j=0}^{d} e^{ik\omega \tau}\, \tau^j \, \bbeta_{k,j},
\end{equation*}
of $\balpha(\tau)$ and $\bbeta(\tau)$ are obtained.

\subsubsection{High precision numerical propagation of a geostationary orbit}

We now consider the initial state of a geostationary satellite \cite[pg. 116]{Montenbruck}.  The initial position coordinates (in kilometers) and velocity vectors (in kilometers per second) are
\begin{align*}
{\bf q}_0 &=(42149.1336,0,0), \\
 \dot{\bf q}_0 &=(0,3.075823259987749,0.0010736649055318406).
\end{align*}
Recall from the discussion above that the  {initial} conditions (\ref{eq:alphabeta_icond})  for the equations (\ref{eq:ode_alphabeta}) can be computed from (\ref{eq:q0}), (\ref{eq:uprime0}) and (\ref{eq:xi}).

%

We have computed several  $(M,d)$-Taylor-Fourier  approximations to the solution of the VOP formulation (\ref{eq:ode_alphabeta}) over 400 revolutions, and found that $M=8$ and $d=8$ is a good choice in that particular initial value problem.
We have estimated the errors of the Taylor-Fourier approximations by comparing {them with} a reference solution computed  (in high precision floating point arithmetic) with Vern9, an explicit RK method of order 9 due to  Verner~\cite{Verner}, efficiently implemented with adaptive time-stepping in DifferentialEquations.jl. The absolute and relative tolerances have been set equal to $10^{-20}$.

\begin{figure}[H]
\centering
{
\includegraphics[width=0.90\hsize]{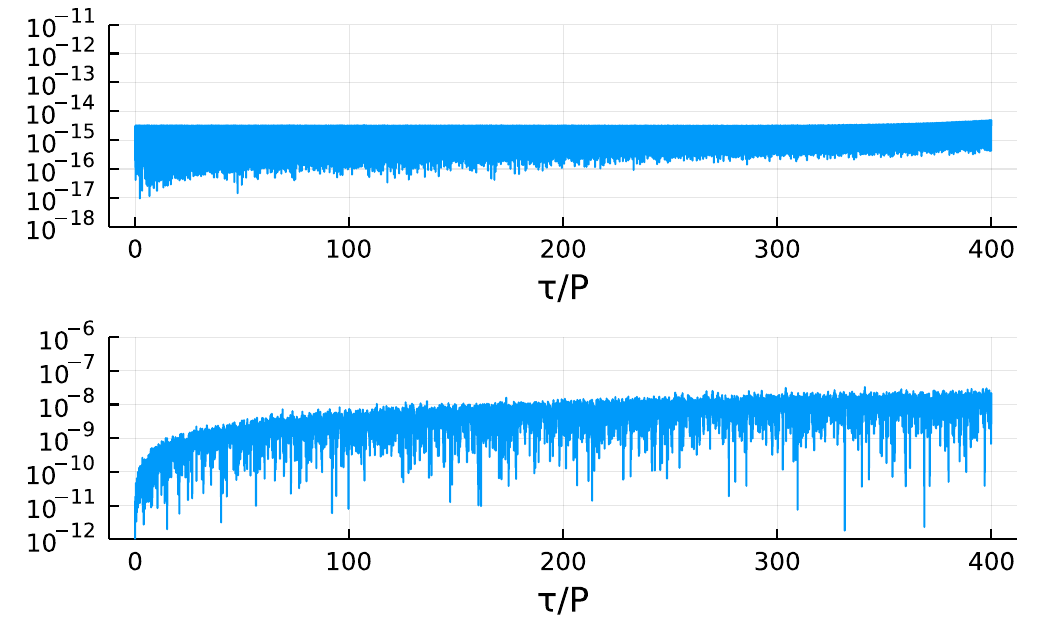}
}
\caption{\label{fig:errors_sat}
    Relative errors in position (top) and absolute errors in physical time (bottom) versus fictitious time $\tau$ (scaled by $P=2\pi/\omega$) of the $(8,8)$-Taylor-Fourier approximation.}
\end{figure}

In Figure~\ref{fig:errors_sat} we display the relative errors in position ${\bf q}$ (top subplot) and the absolute errors in physical time $t$ (bottom subplot) versus the fictitious time $\tau$ scaled by the period $P=2\pi/\omega$ of the satellite orbits for the $(8,8)$-Taylor-Fourier approximation.

We observe that the Taylor-Fourier approximation maintains positional errors below a threshold of approximately $3\times10^{-15}$ (corresponding to the errors due to truncation of Fourier modes) for approximately 380 revolutions. Using $M=16$ (not shown here) instead of $M=8$,  that threshold lowers to round-off levels. Observe that after $380$ revolutions, these errors start to grow. We have checked  {that} error growth is approximately proportional to $t^{d+1}$. The validity of the Taylor-Fourier approximation can be extended to longer time intervals by increasing the degree $d$. For instance, with $d=9$, the errors in position remain below the above-mentioned threshold for up to 500 revolutions.

\begin{figure}[H]
\centering
{
\includegraphics[width=0.90\hsize]{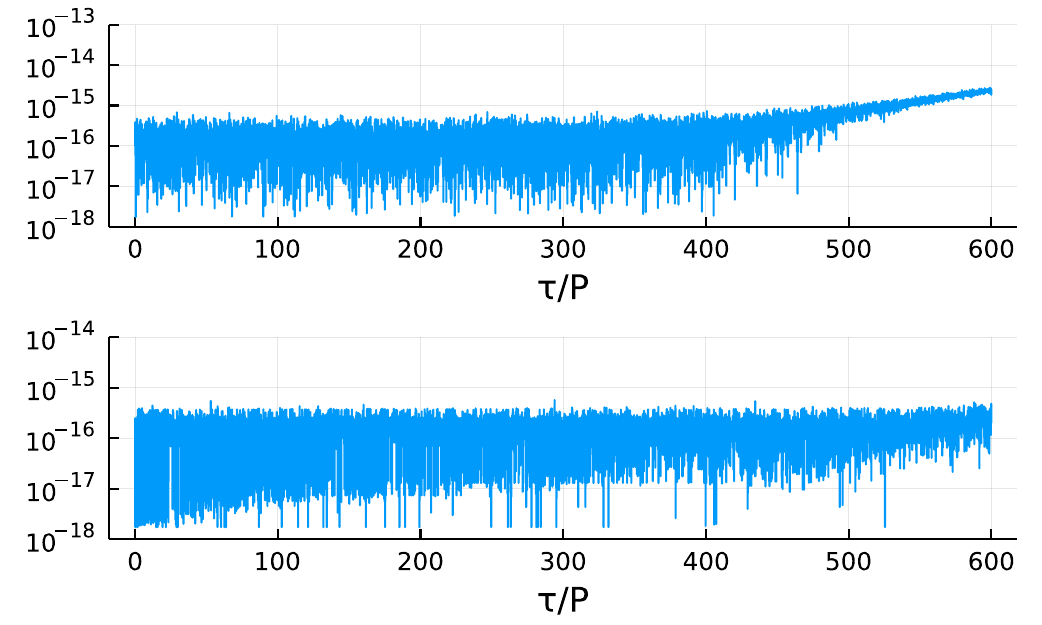}
}
\caption{\label{fig:errors_averaging}
    Energy error of the averaging approximated solution for $(M,d)=(8,8)$ (top) and the position error of the averaged approximation with respect the $(8,8)$-Taylor-Fourier  approximation (bottom) in 600 revolutions.}
\end{figure}

We next check the accuracy of averaged approximation on the left-hand side of
 \eqref{eq:averaging_approximation2}.
We show in Figure~\ref{fig:errors_averaging} the energy error of the averaged approximated solution (top subplot) and the relative error of the position with respect to the Taylor-Fourier approximation (bottom). It can be observed that the behavior of the averaged approximation is practically identical to the behavior of the Taylor-Fourier approximation: it maintains the energy error at rounding error levels for approximately 380 revolutions, just like the Taylor-Fourier approximation, and on the other hand, the difference between the averaged approximation and the Taylor-Fourier approximation remains below the rounding error level for more than 600 revolutions.

\subsubsection{High precision numerical propagation of a highly eccentric orbit}

We now make similar computations with the initial state of a highly eccentric orbit  (with eccentricity $e=0.7679436$) obtained from \cite[pg. 51]{Montenbruck}.
In particular,
\begin{align*}
{\bf q}_0 &=(11959.886901183693,-16289.448826603336,-5963.757695165331), \\
 \dot{\bf q}_0 &=(4.724300951633136,-1.1099935305609756,-0.3847854410416176).
\end{align*}

We have computed several  $(M,d)$-Taylor-Fourier  {approximations to} the solution of the corresponding initial value problem (\ref{eq:ode_alphabeta}), and found that $M=128$ and $d=14$ gives a precision similar to that obtained with $M=8$ and $d=8$ for the geostationary orbit.  This is due to the eccentricity of the  {current} orbit  {which is} much higher than the eccentricity $e \approx 10^{-3}$ of the geostationary orbit,  {and} makes the initial value problem (\ref{eq:ode_alphabeta}) more difficult to  {be solved} numerically (in the sense of requiring more computational effort to achieve similar accuracy). An additional implication of the higher eccentricity is that the interval of validity of the Taylor-Fourier approximations is reduced.

\begin{figure}[H]
\centering
{
\includegraphics[width=0.90\hsize]{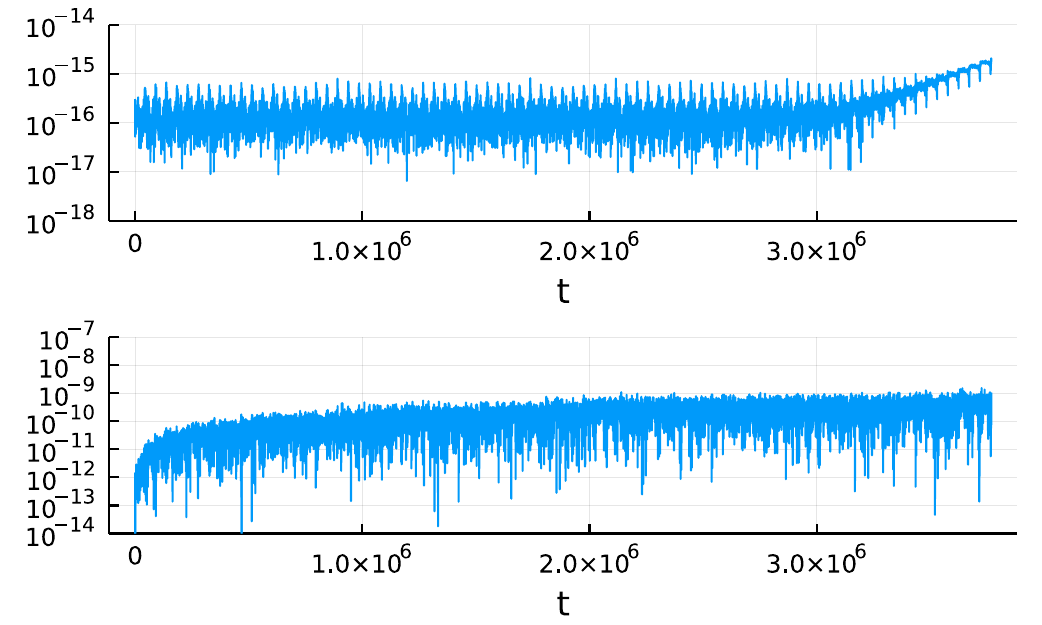}
}
\caption{\label{fig:errors_sat_ecc}
    Relative  {errors in position} (top) and absolute errors in physical time (bottom) versus fictitious time $\tau$ (scaled by $P=2\pi/\omega$) of the $(128,14)$-Taylor-Fourier approximation.}
\end{figure}

Figure~\ref{fig:errors_sat_ecc} is similar to Figure~\ref{fig:errors_sat}, the only differences  being that we display the errors for up to 40 revolutions (instead of 400) and that $(M,d)=(128,14)$  (instead of $(M,d)=(8,8)$).


 In this example, the $(128,14)$-Taylor-Fourier approximation maintains position errors below a threshold of approximately $8\times10^{-16}$ up to $35$ revolutions,  and beyond that, exhibits an error growth that is closely proportional to $t^{d+1}$. We have checked that this threshold is not lowered by increasing $M$, unless a higher precision floating point arithmetic is employed.

\section{Conclusions}
\label{sec:conclusions}

We have proposed approximate solutions to semi-linear ordinary differential equations with highly oscillatory solutions that after an appropriate change of variables can be written as a non-autonomous system with $(2\pi/\omega)$-periodic dependence on $t$. These approximate solutions, that can be written in closed form combining truncated Fourier and Taylor series, are valid   over a time interval whose length is independent of $\omega$. Once computed, they can be efficiently evaluated at arbitrary  times within that interval. The resulting approximate solutions are uniformly accurate~\cite{Ch-C-L-M} with respect to the basic frequency $\omega$.

 We have provided a detailed procedure to compute efficiently such approximate solutions making use of power series arithmetic techniques and the FFT algorithm.  Additionally, we have demonstrated how Taylor-Fourier approximations can be used to evaluate various maps associated with high-order averaging of highly oscillatory systems. Furthermore, we have extended the Taylor-Fourier approximation methodology to highly oscillatory ODEs with multiple non-resonant frequencies.

We have also reported numerical experiments to illustrate the application of the Taylor-Fourier approximations.
Their uniform accuracy with respect to the basic frequency $\omega$ suggests that they could be more efficient than standard numerical methods (which typically suffer from accuracy deterioration as $\omega$ increases) for sufficiently high values of $\omega$.  When applied to a perturbed Kepler problem  (the so-called $J_2$-problem in orbit propagation),  the proposed approximations efficiently provide very accurate solutions, even for highly eccentric orbits. The Taylor-Fourier approximations to the solution of the semi-discretized cubic nonlinear Schr\"odinger equation reproduce the well-known Talbot phenomenon,
though at the cost of increased computational effort due to the high dimensionality of the semi-discretized problem.

\bigskip

{\bf Acknowledgements.} All the authors have received funding by the Spanish State Research Agency through project PID2022-136585NB-C22, \\
 MCIN/AEI/10.13039/501100011033, European Union. M.P.C. was also supported by project VA169P20 (Junta de Castilla y Le\'on, ES) cofinanced by FEDER funds (EU). J.M. and A.M. were partially supported by the Department of Education of the Basque Government through the Consolidated Research Group MATHMODE (ITI456-22).


\begin{thebibliography}{10}

\bibitem{Julia}J. Bezanson, A. Edelman, S. Karpimski \& V.B. Shah, {\em Julia: A Fresh Approach to Numerical computing}, SIAM Review 59, 65--98 (2017).


\bibitem{B-C-M-R2017}  L. Brugnano, M. Calvo, J.I. Montijano \& L. R\'andez, {\em Fourier methods for oscillatory differential problems with a constant high frequency}, AIP Conference Proc. 1863, 020003 (2017). 


\bibitem{C-GP_2015} B. Cano \& A. Gonz\'alez-Pach\'on, {\em Exponential time integration of solitary waves of cubic Schr\"odinger equation}, Appl. Numer. Math. 91, 26--45 (2015)


\bibitem{Ch-M-SS2012a} Ph. Chartier, A. Murua and J. M. Sanz-Serna, {\em Higher-order averaging, formal series and numerical integration II: the quasi-periodic case}, Foundations of Computational Mathematics, 12 (2012), 471-508.

\bibitem{Ch-M-SS2012} P. Chartier, A. Murua \& J.M. Sanz-Serna, {\em A formal series approach to averaging: exponentially small error estimates}, Discrete and Continuous Dynamical Systems 32,  3009--3027 (2012).

\bibitem{Ch-C-L-M} P. Chartier, N. Crouseilles, M. Lemou \& F. M\'ehats, {\em Uniformly accurate numerical schemes for highly oscillatory Klein–-Gordon and nonlinear Schr\"odinger
equations}, Numer. Math. 129,  211--250 (2015).

\bibitem{Chen-Olver}G. Chen \& P.J. Olver, {\em Numerical simulation of nonlinear dispersive quantization}, Discrete and Continuous Dynamical Systems 34, 991--1008 (2014).

\bibitem{C-H-L} D. Cohen, E. Hairer \& Ch. Lubich, {\em Modulated Fourier Expansions of Highly Oscillatory Differential Equations}, Found. Comput. Math. 3, 327--345 (2003).

\bibitem{DP5} J.R. Dormand \& P.J. Prince, {\em A family of embedded Runge-–Kutta formulae}, J. Comput. Appl. Math. 6, 19--26 (1980).

\bibitem{H-L-W} E. Hairer,  Ch. Lubich \& G. Wanner, {\em Geometric Numerical Integration: Structure-Preserving Algorithms for Ordinary Differential Equations}. 2nd ed., Springer, 2006.

\bibitem{H-O2010} M. Hochbruck \& A. Ostermann, {\em Exponential Integrators}, Acta Numerica 19, 209--286 (2010).

\bibitem{K-S1965} P. Kustaanheimo \& E. Stiefel, {\em Perturbation Theory of Kepler Motion Based on Spinor Regularization}, J. Reine Angew. Math. 218, 204-–219 (1965).

\bibitem{Lawson} J.D. Lawson, {\em Generalized Runge–Kutta processes for stable systems with large Lipschitz constants}, SIAM J. Numer. Anal. 4, 372–380 (1967).

\bibitem{Montenbruck} O. Montenbruck \& E. Gill,  {\em Satellite Orbits. Models, Methods and Applications}, Springer, Berlin (2000).

\bibitem{Ostermann-Schratz} A. Ostermann \& K. Schratz, {\em Low regularity exponential-type integrators for semilinear Schr\"odinger equations}, Foundations of Computational Mathematics 18, 731–-755 (2018).

\bibitem{Ostermann-Yao}A. Ostermann \& F. Yao, {\em A fully discrete low-regularity integrator for the nonlinear Schr\"odinger equation}, J. Sci. Comput. 91, 9 (2022).

\bibitem{R-N2017} C. Rackauckas \& Q. Nie, {\em Differentialequations.jl--a performant and feature-rich ecosystem for solving differential equations in julia}, Journal of Open Research Software 5(1):15 (2017).

\bibitem{Stiefel-Scheifele} E. Stiefel \& G. Scheifele, {\em Linear and Regular Celestial Mechanics}, Springer
Verlag, Berlin (1971).

\bibitem{Trefethen} Ll.N. Trefethen, {\em Spectral Methods in Matlab}, SIAM, Philadelphia (2000).

\bibitem{Verner} J. Verner, https://www.sfu.ca/~jverner/.

\end{thebibliography}
\end{document}